\documentclass[11pt]{amsart}
\usepackage[margin=1in]{geometry}
\usepackage{graphicx}
\usepackage{tikz-cd}
\usepackage{amsmath,amssymb,amsthm,mathtools,stmaryrd,enumitem,longtable,booktabs,array,needspace,xcolor,microtype,physics}
\definecolor{citationgreen}{HTML}{07eb57}
\usepackage[colorlinks=false, linkbordercolor=red, citebordercolor=citationgreen, urlbordercolor=cyan,pdfborder={0 0 1}]{hyperref}
\usepackage[nameinlink,capitalise]{cleveref}
\crefformat{equation}{#2(#1)#3}
\Crefformat{equation}{#2(#1)#3}
\numberwithin{equation}{section}

\newtheorem{theorem}{Theorem}[section]
\newtheorem{proposition}[theorem]{Proposition}
\newtheorem{lemma}[theorem]{Lemma}
\newtheorem{corollary}[theorem]{Corollary}
\theoremstyle{definition}
\newtheorem{definition}[theorem]{Definition}

\theoremstyle{remark}

\newcommand{\Z}{\mathbb Z}
\newcommand{\width}{\operatorname{width}}
\newcommand{\Int}{\operatorname{Int}}
\newcommand{\R}{\mathbb R}

\newcommand{\T}{\mathbb T}
\newcommand{\BM}{\mathrm{BM}}
\newcommand{\Sph}{\mathrm{Sph}}

\newcommand{\dist}{\mathrm{dist}}
\newcommand{\udeg}{\mathrm{deg}^{u_X}}
\newcommand{\HH}{\mathcal H}
\newcommand{\cD}{\mathcal{D}}
\newcommand{\cB}{\mathcal B}
\newcommand{\spt}{\mathrm{supp}}
\newcommand{\reg}{\mathrm{reg}}
\newcommand{\sgn}{\mathrm{sgn} }
\title
{The Schoen--Yau--Schick Obstruction in All Dimensions}
\author{Ruotong Jia}
\address{Key Laboratory of Pure and Applied Mathematics, School of Mathematical Sciences, Peking University, Beijing, 100871, P. R. China}
\email{2601110010@stu.pku.edu.cn}
\makeatletter
\renewcommand{\subsection}{%
  \@startsection{subsection}{2}{\z@}%
  {.5\linespacing\@plus.7\linespacing}%
  {.3\linespacing}%
  {\normalfont\bfseries}}
\makeatother
\begin{document}

\begin{abstract}
We prove that closed Schoen--Yau--Schick manifolds admit no metric of positive scalar curvature in all dimensions. Building on the singular dimension descent of Bi--Zhu, we extend their topological argument from overtorical manifolds to the Schoen--Yau--Schick setting. The key is to detect a non-spherical two-dimensional homology class by a cohomology class with finite coefficients and preserve the resulting nonzero degree throughout the descent, leading to a contradiction on the terminal surface.
\end{abstract}

\maketitle

\tableofcontents

\section{Introduction}\label{sec:introduction}
By the Gauss--Bonnet theorem, a closed connected oriented surface admitting a metric of positive Gauss curvature is diffeomorphic to the sphere. The minimal hypersurface method of Schoen--Yau provides a way \cite{SY79} to exploit this restriction in higher dimensions: successive codimension--one reductions connect scalar curvature to the topology of a terminal surface, which gives a proof of Geroch conjecture in dimensions $n\le 7$, i.e. the torus $\T^n$ admits no metric of positive scalar curvature.

The Schoen--Yau--Schick (SYS) condition identifies a homological obstruction suited to this procedure. In dimensions $3\le n\le 7$, the SYS obstruction follows from the minimal hypersurface method of Schoen--Yau \cite{SY79}, by iterating Schick's cap product theorem \cite[Corollary 1.5]{Sch98}.

In higher dimensions, minimizing hypersurfaces may be singular, and continuing the descent requires both geometric control near the singular sets and a means of retaining the relevant homological information. Schoen--Yau \cite{SY17} proposed an extension of the minimal hypersurface method to all dimensions using minimal slicings.

Recent work of Bi--Zhu \cite{BZ26} has shown that closed oriented enlargeable manifolds admit no metric of positive scalar curvature. Their argument proceeds by dimension descent in the presence of singularities, with the dimension bounds for the singular sets decreasing at each step. By combining weighted scalar curvature estimates with a carefully controlled conformal blow--up near the singular sets, they retain the curvature control needed to continue the descent. Bi--Zhu also discuss dimension descent with an integral cohomological degree in Remarks 4.12 and 5.2. For SYS manifolds, we use a finite-coefficient cohomology class to detect the non-spherical two-dimensional homology class and show that the resulting nonzero pairing is preserved at each step of the descent. Our main result is the following.
\begin{theorem}\label{thm:Main}
Let $n\ge3$ and let $M^n$ be a closed SYS manifold. Then $M$ admits no metric of positive scalar curvature.
\end{theorem}

The descent applies more generally to complete manifolds carrying a proper map of nonzero detected degree.
\begin{theorem}\label{thm:abstract-obstruction}
Let $(X^n,g_X)$ be a connected complete oriented Riemannian manifold without boundary, $n\ge2$, with $R_{g_X}\ge\sigma>0$. Let $A=\Z/q$, $q\ge2$, and let $u_X\in H^2(X;A)$ satisfy
$$
\langle u_X,w\rangle=0,\ \forall w\in\Sph_2(X).
$$
For every smooth proper map $\Phi:X\to\R^{n-2}$ whose coordinates are globally Lipschitz, one has $\udeg(\Phi)=0$.
\end{theorem}
Here $\Sph_2(X)$ is the image of the Hurewicz map, and the degree is defined by
$$
\Phi_*\left([X]_{\BM,\Z}\frown u_X\right)=\udeg(\Phi)[\R^{n-2}]_{\BM,A}.
$$
For a regular value $z$, it equals $\langle u_X|_{\Phi^{-1}(z)},[\Phi^{-1}(z)]_\Z\rangle$, see \cref{def:u_X degree,prop:fibre degree}.

\begin{proof}[Proof of \cref{thm:Main} from \cref{thm:abstract-obstruction}]
Suppose that a closed SYS manifold $(M,g)$ has $R_g>0$. Compactness implies $R_g\ge\sigma>0$. Choose a SYS class $z$ as in \cref{def:SYS} and a detector $u$ for $z$ as in \cref{prop:detector}. By \cref{prop:carrier}, there are a connected covering $p:X\to M$ and a smooth proper map $\Phi:X\to\R^{n-2}$ with globally Lipschitz coordinates and $\udeg(\Phi)\ne0$, where $g_X=p^*g$ and $u_X=p^*u$. The metric $g_X$ is complete and $R_{g_X}\ge\sigma$. Moreover, for every $v:S^2\to X$,
$$
\langle u_X,v_*[S^2]_\Z\rangle=\langle u,(p\circ v)_*[S^2]_\Z\rangle=0.
$$
This contradicts \cref{thm:abstract-obstruction}.
\end{proof}

It remains to prove \cref{thm:abstract-obstruction}. Assume that $\udeg(\Phi)\ne0$. The construction below starts directly from $(X,g_X,\Phi,u_X)$ and preserves this degree throughout the descent.

At each stage, we construct a wide band and a weighted $\mu$-bubble separating its two faces. The topological step compares the two-dimensional fibres of the bubble with those of a reference hypersurface. After choosing fibres which avoid the singular sets, slicing a comparison current shows that these surfaces represent the same homology class in the incoming regular manifold. Their pairings with $u_X$ therefore agree, preserving the nonzero degree at the next stage.

For the geometric step, we follow the weighted stability and conformal blow-up method of Bi--Zhu. We construct a complete metric on the regular part of the bubble with a positive weighted scalar curvature lower bound, while the dimension bound for the singular set decreases by one. Iterating the construction reaches a closed smooth surface. The curvature inequality and Gauss--Bonnet force every component to be a sphere, on which the detector vanishes, contradicting the preserved nonzero degree.

\section{Preliminaries}\label{sec:preliminaries}
\subsection{SYS manifolds and the detector}\label{subsec:sysmanifolds}
\begin{definition}[Spherical Classes]\label{def:spherical classes}
For a connected manifold $M$ and an integer $k\ge1$, define the subgroup of spherical $k$-classes by
$$
\Sph_k(M)=\mathrm{im}\left(\pi_k(M)\longrightarrow H_k(M;\Z)\right).
$$
\end{definition}
\begin{proposition}[The detector for non-spherical classes]\label{prop:detector}
Let $M$ be a closed connected manifold, $k\ge1$, and $z\in H_k(M;\Z)\setminus\Sph_k(M)$. Then there are $q\ge2$ and $u\in H^k(M;A)$, with $A=\Z/q$, such that
$$
\langle u,z\rangle\ne0,\ \langle u, w\rangle=0,\ \forall w\in\Sph_k(M).
$$
\end{proposition}
\begin{proof}
The group $G=H_k(M;\Z)/\Sph_k(M)$ is finitely generated abelian. Write
$$
G\cong\Z^r\oplus\Z/d_1\oplus\cdots\oplus\Z/d_j.
$$
The image of $z$ has a nonzero coordinate. If a free coordinate is the nonzero integer $a_0$, reduce that coordinate modulo an integer $q=|a_0|+1$. If a torsion coordinate is nonzero in $\Z/d_i$, project to that coordinate and take $q=d_i$. In either case there is a homomorphism $\chi:G\to A$ which does not vanish on the image of $z$.

Compose it with projection $H_k(M;\Z)\to G$, obtaining $\ell:H_k(M;\Z)\to A$. By the universal coefficient theorem, the evaluation map
$$
H^k(M; A)\longrightarrow \mathrm{Hom}_{\Z}\left(H_k(M;\Z),A\right)
$$
is surjective. Choose $u$ mapping to $\ell$. Then $\langle u,z\rangle=\ell(z)\ne0$, whereas $\langle u,w\rangle=\ell(w)=0$ for every $w\in\Sph_k(M)$.
\end{proof}
The SYS condition specifies a two-dimensional integral homology class outside the spherical subgroup.
\begin{definition}[SYS manifolds]\label{def:SYS}
Let $n\ge 3$, a closed connected oriented manifold $M^n$ is called SYS (Schoen--Yau--Schick) if there are $\alpha_1,\cdots,\alpha_{n-2}\in H^1(M;\Z)$ with
$$
z=[M]_{\Z}\frown\left(\alpha_1\smile\cdots\smile\alpha_{n-2}\right) \notin\Sph_2(M).
$$
\end{definition}
For the class $z$ in \cref{def:SYS}, choose $u\in H^2(M;A)$ as in \cref{prop:detector}, so that $\langle u,z\rangle\ne0$ and $\langle u,w\rangle=0$ for every $w\in\Sph_2(M)$. We call $u$ a detector for $z$.

To see the relation between SYS manifolds and the torus, one can construct a smooth map $\Theta=(\Theta_1, \cdots, \Theta_{n-2}): M\to \T^{n-2}$ as follows:

Identify $S^1$ with $\R/\Z$, and let $d\theta$ denote the global $1$-form so that $\int_{S^1}d\theta=1$. Let $\eta\in H^1(S^1;\Z)$ be the positive generator, whose image in real cohomology is represented by $d\theta$. For each $i=1,\cdots,n-2$, choose a smooth map $\Theta_i:M\to S^1$ such that $\Theta_i^*\eta=\alpha_i$.

Now the SYS condition should be a topological condition on $\Theta$, which is expressed by so called $u$--degree. Consider $[M]_{\Z}\frown u\in H_{n-2}(M; A)$, there exists a constant in $A$, denoted by $\deg^u(\Theta)$ such that
$$\Theta_*\left([M]_{\Z}\frown u\right)=\deg^u(\Theta)[\T^{n-2}]_A.$$
\begin{proposition}\label{prop:degcompute}
$\deg^u(\Theta)=\langle u,z\rangle\ne 0$.
\end{proposition}
\begin{proof}
For $i=1, \cdots, n-2$, write $\pi_i: \T^{n-2}\to S^1$ to be the projection to the $i$--th coordinate, define $\eta_{\T^{n-2}}=\pi_1^*\eta\smile\cdots\smile\pi_{n-2}^*\eta\in H^{n-2}(\T^{n-2};\Z)$, then $\langle \eta_{\T^{n-2}}, [\T^{n-2}]_{\Z}\rangle=1$, and
$$
\begin{aligned}
\deg^u(\Theta) &=\bigl\langle\eta_{\T^{n-2}}, \Theta_*\left([M]_{\Z}\frown u\right)\bigr\rangle\\
&=\bigl\langle\Theta^*\eta_{\T^{n-2}},[M]_{\Z}\frown u\bigr\rangle\\
&=\bigl\langle u\smile\Theta^*\eta_{\T^{n-2}},[M]_{\Z}\bigr\rangle\\
&=\bigl\langle u\smile\alpha_1\smile\cdots\smile\alpha_{n-2}, [M]_{\Z}\bigr\rangle\\
&=\bigl\langle u,[M]_{\Z}\frown (\alpha_1\smile\cdots\smile\alpha_{n-2})\bigr\rangle\\
&=\langle u,z\rangle\ne0.
\end{aligned}
$$
\end{proof}
\subsection{Borel--Moore fundamental classes and the degree formula}\label{subsec:BM}
Motivated by \cref{prop:degcompute}, we extend the notion of $u$-degree to proper maps between possibly noncompact manifolds. We first recall the locally finite singular-chain model of Borel--Moore homology below. All manifolds in this subsection are without boundary.

Let $R$ be a commutative ring with unity and $M$ be an $n$-manifold. A locally finite singular $k$-chain with coefficients in $R$ is a possibly infinite formal sum
$$
c=\sum_{\lambda}a_{\lambda}\sigma_{\lambda},
$$
where $a_\lambda\in R$ and $\sigma_\lambda:\Delta^k\to M$ are singular simplices, such that every compact subset of $M$ meets the images of only finitely many simplices with nonzero coefficient. The group of such chains is denoted by $C_k^{\BM}(M;R)$. The usual boundary operator preserves local finiteness and satisfies $\partial^2=0$, so these groups form a chain complex $C_*^{\BM}(M;R)$. Its homology is the Borel--Moore homology $H_*^{\BM}(M;R)$.

If further $M$ is oriented, for $x\in M$, choose the generator $o_x$ of the $R$-module $H_n\left(M,M\setminus\left\{x\right\};R\right)$ determined by the orientation of $M$. For $c\in C_k^{\BM}(M;R)$ as above, write
$$
L_x(c)=\left[\sum_{\lambda:x\in\sigma_\lambda(\Delta^k)} a_\lambda\sigma_\lambda\right] \in C_k\left(M,M\setminus\left\{x\right\};R\right),
$$
where the brackets denote the coset in $C_k(M;R)/C_k\left(M\setminus\left\{x\right\};R\right)$. The sum is finite by local finiteness. A simplex missing $x$, together with all its faces, vanishes in this quotient, hence $L_x\partial=\partial L_x$ and $L_x$ induces a homomorphism
$$
\rho_x: H_n^{\BM}(M;R)\longrightarrow H_n\left(M,M\setminus\left\{x\right\};R\right)
$$

\begin{proposition}[The BM fundamental class]\label{prop:BM-fundamental}
An oriented $M$ has a unique class $[M]_{\BM,R}\in H_n^{\BM}(M;R)$ such that
$$
\rho_x\left([M]_{\BM,R}\right)=o_x\ \text{for every }x\in M.
$$
\end{proposition}
\begin{proof}
Take $X=A=M$ in \cite[Proposition 2.12(3)]{CZ26}.
\end{proof}

If $M$ is closed, compactness forces every locally finite chain to be finite. Thus $H_*^{\BM}(M;R)=H_*(M;R)$, and $[M]_{\BM,R}=[M]_{R}$. On a noncompact oriented $n$-manifold, the BM class is the appropriate replacement of the fundamental class.

\begin{definition}[The $u_X$--degree]\label{def:u_X degree}
Fix $A=\Z/q$, $q\ge2$. Let $X$ be a topological space with a class $u_X\in H^2(X;A)$, $N^{d+2}$ be an oriented manifold with a continuous map $\iota:N\to X$, $V^d$ be a connected oriented manifold and $F:N\to V$ be a proper continuous map. Define
\begin{equation}\label{equ:u_X degree}
\udeg(F)[V]_{\BM,A} =F_*\left([N]_{\BM,\Z}\frown\iota^*u_X\right)\in H_d^{\BM}(V;A).
\end{equation}
\end{definition}
\begin{lemma}
\label{lem:degree-surjective}
If $\udeg(F)\ne0$, then $F$ is surjective.
\end{lemma}
\begin{proof}
Suppose $F$ is not surjective. Properness makes $F(N)$ closed, so there is a coordinate ball $B\subset V\setminus F(N)$. Fix $y\in B$ and choose a compactly supported cocycle representing $\tau_y\in H_c^d(V; A)$, with support contained in $B$, normalized by
$$
[V]_{\BM, A}\frown\tau_y=[y]_{A}.
$$
Then
$$
\begin{aligned}
\udeg(F)[y]_A&=\udeg(F)[V]_{\BM, A}\frown\tau_y\\
&=\left(F_*\left([N]_{\BM,\Z} \frown\iota^*u_X\right)\right)\frown\tau_y\\
&=F_*\left(\left([N]_{\BM,\Z} \frown\iota^*u_X\right)\frown F^*\tau_y\right)\\
&=0,
\end{aligned}
$$
a contradiction.
\end{proof}
\begin{proposition}[The degree formula]\label{prop:fibre degree}Assume that $F$ is smooth. Then for every regular value $y\in V$, its closed oriented fibre $\Gamma_y=F^{-1}(y)$ satisfies
$$
\udeg(F)=\langle (\iota|_{\Gamma_y})^*u_X,[\Gamma_y]_{\Z}\rangle.
$$
\end{proposition}
\begin{proof}
Fix a regular value $y\in V$. If $\Gamma_y$ is empty, both sides vanish by \cref{lem:degree-surjective} and the convention $[\varnothing]_\Z=0$. Otherwise, orient $\Gamma_y$ so that a positive frame tangent to the fibre followed by lifts of a positive frame on $V$ is positive in $N$. Since $F$ is proper, its critical values form a closed set. Choose a coordinate ball $B$ about $y$ containing no critical values. The restriction $F_B:F^{-1}(B)\to B$ is a proper submersion whose image is nonempty, open and closed in $B$. By Ehresmann's fibration theorem \cite[Lemma 9.2]{KSM93}, after shrinking $B$ if necessary, choose a smooth trivialization
$$
\Psi:F^{-1}(B)\to\Gamma_y\times B.
$$
Let $p_\Gamma:\Gamma_y\times B\to\Gamma_y$ and $p_B:\Gamma_y\times B\to B$ be the projections, and let $j:\Gamma_y\to\Gamma_y\times B$ be the inclusion $x\mapsto (x, y)$. The trivialization satisfies $p_B\circ\Psi=F_B$ and $j=\Psi|_{\Gamma_y}$.

Since $B$ is contractible, choose a homotopy $h: B\times[0, 1]\to B$, such that
$$h(\cdot, 0)=Id_B,\ h(\cdot, 1)=y.$$

Then $((x, b), t)\mapsto (x, h(b, t))$ is a homotopy between $Id_{\Gamma_y\times B}$ and $j\circ p_\Gamma$.

Consequently,
$$
p_\Gamma^*(\iota|_{\Gamma_y})^*u_X=p_\Gamma^*(\Psi^{-1}\circ j)^*(\iota|_{F^{-1}(B)})^*u_X=(\Psi^{-1})^*(\iota|_{F^{-1}(B)})^*u_X.
$$
The construction of the BM fundamental class implies
$$
\Psi_*[F^{-1}(B)]_{\BM,\Z}=[\Gamma_y\times B]_{\BM, \Z} =[\Gamma_y]_{\Z}\times[B]_{\BM,\Z}.
$$
Let $\tau_B\in H_c^d(B;\Z)$ be the orientation class normalized by $[B]_{\BM,\Z}\frown\tau_B=[y]_\Z$. Its extension by zero $\tau_y\in H_c^d(V;\Z)$ satisfies
$$
[V]_{\BM, \Z}\frown\tau_y=[y]_{\Z}.
$$
Then we obtain
$$
F_*\left(\left([N]_{\BM,\Z}\frown\iota^*u_X\right) \frown F^*\tau_y\right)=\udeg(F)[y]_A.
$$
Use the inclusion $B\hookrightarrow V$ to identify $H_0(B;A)$ with $H_0(V;A)$. Since $\tau_y$ is extended from $B$, excision for the compactly supported cap product yields
$$
\begin{aligned}
&F_*\left(\left([N]_{\BM,\Z}\frown\iota^*u_X\right)\frown F^*\tau_y\right)\\
&\qquad=(F_B)_*\left(\left([F^{-1}(B)]_{\BM,\Z}\frown(\iota|_{F^{-1}(B)})^*u_X\right)\frown F_B^*\tau_B\right).
\end{aligned}
$$

Naturality of the cap product, applied twice, now yields
$$
\begin{aligned}
&\Psi_*\left( \left([F^{-1}(B)]_{\BM,\Z}\frown (\iota|_{F^{-1}(B)})^*u_X\right) \frown F_B^*\tau_B\right)\\
& =\Psi_*\left( \left([F^{-1}(B)]_{\BM,\Z}\frown \Psi^*p_\Gamma^*(\iota|_{\Gamma_y})^*u_X\right) \frown\Psi^*p_B^*\tau_B\right)\\
& =\Psi_*\left( [F^{-1}(B)]_{\BM,\Z}\frown \Psi^*p_\Gamma^*(\iota|_{\Gamma_y})^*u_X\right) \frown p_B^*\tau_B\\
& =\left( \Psi_*[F^{-1}(B)]_{\BM,\Z}\frown p_\Gamma^*(\iota|_{\Gamma_y})^*u_X\right) \frown p_B^*\tau_B\\
& =\left( \left([\Gamma_y]_{\Z}\times[B]_{\BM,\Z}\right) \frown p_\Gamma^*(\iota|_{\Gamma_y})^*u_X\right) \frown p_B^*\tau_B.
\end{aligned}
$$
The cap--cross product formula separates the detector from the base factor:
$$
\left([\Gamma_y]_{\Z}\times[B]_{\BM,\Z}\right) \frown p_\Gamma^*(\iota|_{\Gamma_y})^*u_X = \left([\Gamma_y]_{\Z}\frown(\iota|_{\Gamma_y})^*u_X\right) \times[B]_{\BM,\Z}.
$$
Therefore,
$$
\begin{aligned}
\udeg(F)[y]_A &=(p_B)_*\left(\left(([\Gamma_y]_{\Z}\times[B]_{\BM,\Z}) \frown p_\Gamma^*(\iota|_{\Gamma_y})^*u_X\right) \frown p_B^*\tau_B\right)\\
&=(p_B)_*\left(\left( ([\Gamma_y]_{\Z}\frown(\iota|_{\Gamma_y})^*u_X) \times[B]_{\BM,\Z}\right)\frown p_B^*\tau_B\right)\\
&=(p_B)_*\left( \left([\Gamma_y]_{\Z}\frown(\iota|_{\Gamma_y})^*u_X\right) \times\left([B]_{\BM,\Z}\frown\tau_B\right)\right)\\
&=(p_B)_*\left( \left([\Gamma_y]_{\Z}\frown(\iota|_{\Gamma_y})^*u_X\right) \times[y]_{\Z}\right)\\
&=\langle(\iota|_{\Gamma_y})^*u_X,[\Gamma_y]_{\Z}\rangle[y]_A.
\end{aligned}
$$
Comparing with $\udeg(F)[y]_A$ and using $H_0(V;A)\cong A$ proves the assertion.
\end{proof}
\subsection{The initial covering construction}
\label{sec:carrier}
\begin{proposition}[A connected proper carrier]\label{prop:carrier}
Let $n\ge3$, $(M^n, g)$ be an SYS manifold, and fix a SYS class $z$ and a detector $u\in H^2(M;A)$ for $z$ as in \cref{def:SYS,prop:detector}. There exists a connected covering $p:X\to M$, the lifted complete metric $g_X=p^*g$, and a smooth proper surjective map
$$
\Phi=(\Phi_1,\cdots,\Phi_{n-2}):X\to\R^{n-2},
$$
whose coordinates are globally Lipschitz and for which
$$
u_X=p^*u,\ \udeg(\Phi)\ne0\in A.
$$
\end{proposition}
$$
\begin{tikzcd}[row sep=large, column sep=large]
X \arrow[r, "\Phi"] \arrow[d, "p"'] & \R\mathrlap{{}^{n-2}} \arrow[d, "\pi_\T"] \\
M \arrow[r, "\Theta"'] & \T\mathrlap{{}^{n-2}}
\end{tikzcd}
$$
\begin{proof}
Denote $\pi_\T$ by the universal covering
$$\pi_\T: \R^{n-2}\to \T^{n-2}$$ 
For $i=1,\cdots,n-2$, let $d\theta_i=\pi_i^*d\theta$ be the standard closed $1$-form on $\T^{n-2}$ associated with the $i$-th circle factor. Define
$$
\widehat X=\left\{(x,y)\in M\times\R^{n-2}:\Theta(x)=\pi_{\T}(y)\right\},\ \widehat p(x,y)=x,\ \widehat\Phi(x,y)=y, g_{\widehat X}=\widehat p^*g.
$$
This is a smooth covering of $M$ and $\widehat\Phi$ is proper. The lifted metric is complete componentwise. Moreover
$$
\sup_{\widehat X}|d\widehat\Phi_i|_{g_{\widehat X}}=\sup_{\widehat X}|\widehat p^*\Theta^*d\theta_i|_{g_{\widehat X}}=\sup_M|\Theta^*d\theta_i|_g<\infty,
$$
Integrating along curves yields the global Lipschitz bound on each component.

With the integral orientation class $\eta_{\T^{n-2}}$ used in \cref{prop:degcompute}, the identity
$$
\langle u\smile\Theta^*\eta_{\T^{n-2}},[M]_\Z\rangle=\langle u,z\rangle\ne0
$$
implies $\Theta^*\eta_{\T^{n-2}}\ne0$. If $\Theta$ omitted a point of $\T^{n-2}$, it would factor through the punctured torus, whose ordinary top cohomology in degree $n-2$ is zero. Thus $\Theta$ and $\widehat\Phi$ are surjective.

Choose a regular value $y$ of $\widehat\Phi$. The covering projection $\widehat p$ identifies its oriented fibre with the fibre of $\Theta$ over $\pi_\T(y)$. The latter represents $z$ in $M$, so \cref{prop:fibre degree} yields
$$
\begin{aligned}
\deg^{\widehat p^*u}(\widehat\Phi)
&=\bigl\langle (\widehat p|_{\widehat\Phi^{-1}(y)})^*u,[\widehat\Phi^{-1}(y)]_\Z\bigr\rangle\\
&=\bigl\langle u|_{\Theta^{-1}(\pi_\T(y))},[\Theta^{-1}(\pi_\T(y))]_\Z\bigr\rangle=\langle u,z\rangle.
\end{aligned}
$$

Choose a bounded open ball $B\subset\R^{n-2}$, then $\widehat{\Phi}^{-1}(\overline{B})$ is compact, which meets only finitely many of connected components of $\widehat{X}$, say $X_1, \cdots, X_r$.

The preceding identity implies
$$\begin{aligned}
\deg^{\widehat p^*u}(\widehat{\Phi})[B]_{\BM, A} &=\widehat{\Phi}_*\left([\widehat{X}]_{\BM, \Z}\frown \widehat p^*u\right)|_B \\
&=\sum_{i=1}^r(\widehat{\Phi}|_{X_i})_*\left([X_i]_{\BM, \Z}\frown(\widehat p|_{X_i})^*u\right)|_B   \\
&=\sum_{i=1}^r\deg^{(\widehat p|_{X_i})^*u}(\widehat{\Phi}|_{X_i})[B]_{\BM, A}.
\end{aligned}
$$

Now choose one component $X_i$ with $\deg^{(\widehat p|_{X_i})^*u}(\widehat{\Phi}|_{X_i})\ne 0$, define
$$X=X_i,\ p=\widehat{p}|_X,\ \Phi=\widehat{\Phi}|_X,\ g_X=p^*g,\ u_X=p^*u.$$
\end{proof}
\section{The descent framework}\label{sec:the descent framework}
\subsection{The descent datum}\label{subsec:the descent datum}
Fix a complete oriented Riemannian manifold $(X^n,g_X)$ without boundary, a class $u_X\in H^2(X;A)$ annihilating spherical classes, and a smooth map $\Phi:X\to\R^k$ with globally Lipschitz coordinates. Here $A=\Z/q$ with $q\ge2$. Denote by $d_X$ the distance induced by $g_X$, and by $\HH$ the Hausdorff measure of the metric space under consideration.
\begin{definition}[The $m$-dimensional datum]\label{def:datum}
For $2\le m\le n$ with $m-2\le k$, write
$$
\cD_m=(X_m^\circ,\overline{X_m},D_m,F_m;g_m,f_m,\sigma_m).
$$
The required conditions are:
\begin{enumerate}[label=\textup{(D\arabic*)}]
\item\label{item:strata} $X_m^\circ$ is an oriented smoothly embedded $m$-manifold in $X$. Its closure $\overline{X_m}=\overline{X_m^\circ}^X$. The complement
$$
D_m=\overline{X_m}\setminus X_m^\circ
$$
is closed, and $\dim_\HH D_m\le m-7$. Here, whenever a Hausdorff dimension bound is negative, the corresponding set is understood to be empty.
\item\label{item:map} The map
$$
F_m=(\Phi_1,\cdots,\Phi_{m-2})|_{\overline{X_m}}: \overline{X_m}\longrightarrow\R^{m-2}
$$
is proper. It is locally Lipschitz with respect to $d_X$, and smooth on $X_m^\circ$.
\item\label{item:degree} Write
$$
V_m=\R^{m-2}\setminus F_m(D_m),\ N_m=F_m^{-1}(V_m)\subset X_m^\circ,
$$
and regard $F_m|_{N_m}$ as a proper map $N_m\to V_m$. With $\iota:N_m\hookrightarrow X$ the inclusion, require
$$
\udeg(F_m|_{N_m})\ne0.
$$
Note that the dimension bound $\dim_{\mathcal H}F_m(D_m)\le m-7$ implies that $V_m$ is connected. 
\item\label{item:metric} $g_m$ is a smooth complete metric on every
component of $X_m^\circ$, dominating the metric induced from $g_X$. The function $f_m$ is smooth, and
$$
\mathcal S_{f_m}^{m,-1}(g_m) =R_{g_m}+2\Delta_{g_m}f_m-\frac{m}{m+1}|df_m|_{g_m}^2 \ge\sigma_m>0.
$$
\item\label{item:inverse} For every $a\in D_m$, there are
$r_a,c_a>0$ such that for $x\in X_m^\circ\cap B_{d_X}(a,r_a)$ and $v\in T_xX_m^\circ$,
$$
g_m(v,v)\ge\frac{c_a}{d_X(x,a)^2}g_X(v,v).
$$
\end{enumerate}
\end{definition}
\subsection{The reference slice}\label{subsec:the reference slice}
Fix an $m$-dimensional datum $\cD_m$ with $m\ge3$. Separate the last coordinate by writing
$$
F_m=(G_m,t_m),\ G_m=(\Phi_1,\cdots,\Phi_{m-3})|_{\overline{X_m}},\ t_m=\Phi_{m-2}|_{\overline{X_m}}.
$$
\begin{lemma}\label{lem:goodlevel}
One can choose $c_m\in\R$ so that
$$
P_m=X_m^\circ \cap t_m^{-1}(c_m)=\left(t_m|_{X_m^\circ}\right)^{-1}(c_m)
$$
is a smooth oriented hypersurface and
$$
\dim_\HH\left(D_m\cap t_m^{-1}(c_m)\right)\le m-8.
$$
Moreover, for almost every $z\in\R^{m-3}$, $(z,c_m)\notin F_m(D_m)$ is a regular value of $F_m|_{X_m^\circ}$ and
$$
\langle u_X|_{P_{m,z}},[P_{m,z}]_\Z\rangle=\udeg(F_m|_{N_m}),
$$
where $P_{m,z}=P_m\cap G_m^{-1}(z)=X_m^\circ\cap F_m^{-1}(z,c_m)$.
\end{lemma}
\begin{proof}
Choose a countable open cover $\left\{U_j\right\}_{j\ge1}$ of $\overline{X_m}$ such that $t_m|_{U_j}$ is Lipschitz for any $j\ge 1$.

If $m\le7$, local Lipschitz continuity and countable stability of Hausdorff dimension imply $\dim_\HH t_m(D_m)\le0$. Thus $D_m\cap t_m^{-1}(c)=\varnothing$ for almost every $c\in\R$.

Suppose now that $m\ge8$. For every rational $s>m-7$ and $j\ge 1$, $\HH^s(D_m\cap U_j)=0$ by \ref{item:strata}. The coarea inequality \cite[Theorem 1.1]{EH21}, applied to $t_m|_{U_j}$ with target $\R$, implies
$$
\HH^{s-1}(D_m\cap U_j\cap t_m^{-1}(c))=0
$$
for almost every $c\in\R$. Taking the union of these exceptional sets over the countably many $j$ and rational $s>m-7$, we obtain, for almost every $c$,
$$
\HH^{s-1}\left(D_m\cap t_m^{-1}(c)\right)=0,\ \text{for every rational }s>m-7.
$$
Letting $s$ decrease to $m-7$ proves
$$
\dim_\HH\left(D_m\cap t_m^{-1}(c)\right)\le m-8.
$$

By Sard's theorem, the critical values of $t_m|_{X_m^\circ}$ also form a Lebesgue null set. Choose $c_m\in\R$ outside this set and the preceding exceptional sets. Then $P_m$ is a smooth hypersurface.

Local Lipschitz continuity of $G_m$ and the slicing bound imply
$$
\dim_\HH G_m\left(D_m\cap t_m^{-1}(c_m)\right) \le m-8<m-3.
$$
Consequently, almost every $z\in\R^{m-3}$ satisfies $(z,c_m)\notin F_m(D_m)$. Sard's theorem applied to $G_m|_{P_m}$ further ensures that almost every such $z$ is its regular value, and $(z,c_m)$ is a regular value of $F_m|_{X_m^\circ}$.

For each of these values, $P_{m,z}=F_m^{-1}(z,c_m)$ is a closed oriented surface contained in $N_m$. Applying \cref{prop:fibre degree} to $F_m|_{N_m}$ and using \ref{item:degree}, we conclude that
$$
\langle u_X|_{P_{m,z}},[P_{m,z}]_\Z\rangle =\udeg(F_m|_{N_m})\ne0.
$$
\end{proof}
\begin{proposition}\label{prop:localization}
Let $\dist_{g_m}$ denote the Riemannian distance induced by $g_m$ on $X_m^\circ$, with distance $\infty$ between distinct connected components. For every $W>0$:
\begin{enumerate}[label=\textup{(\roman*)}]
\item There are points $x_-,x_+$ in the same connected component of $X_m^\circ$ such that
$$
t_m(x_-)<c_m<t_m(x_+),\ \dist_{g_m}(x_\pm,P_m)>W.
$$
\item The distance neighborhood
$$
N_W(P_m)=\left\{x\in X_m^\circ:\dist_{g_m}(x,P_m)<W\right\}
$$
satisfies
$$
\overline{N_W(P_m)}^{X}\cap D_m \subset D_m\cap t_m^{-1}(c_m).
$$
In particular, 
$$\dim_\HH\left(\overline{N_W(P_m)}^{X}\cap D_m\right)\le m-8.$$
\end{enumerate}
\end{proposition}
\begin{proof}
Choose a connected component $C_m$ of $N_m$ such that $\udeg(F_m|_{C_m})\ne0$. By \cref{lem:degree-surjective}, $F_m|_{C_m}$ maps onto $V_m$. Since $\dim_\HH F_m(D_m)\le m-7$, we obtain $t_m(C_m)=\R$.

\ref{item:metric} yields
$$\sup_{X_m^\circ}|dt_m|_{g_m}=\sup_{X_m^\circ}|d\Phi_{m-2}|_{g_m}\le\sup_{X}|d\Phi_{m-2}|_{g_X}<\infty.$$
Thus,
$$
\dist_{g_m}(x,P_m)\ge\frac{|t_m(x)-c_m|}{\sup\limits_X|d\Phi_{m-2}|_{g_X}},\ \forall x\in X_m^\circ.
$$
Choose $x_-,x_+\in C_m$ with $t_m(x_-)<c_m-\sup_{X}|d\Phi_{m-2}|_{g_X}W$ and $t_m(x_+)>c_m+\sup_{X}|d\Phi_{m-2}|_{g_X}W$. This proves (i).

For (ii), let $\left\{x_j\right\}\subset N_W(P_m)$ converge in $X$ to $a\in D_m$. Choose $p_j\in P_m$ and paths from $x_j$ to $p_j$ of $g_m$-length at most $W+1$. We claim that $p_j\to a$. Otherwise, after passing to a subsequence, $d_X(p_j,a)\ge r$ for some $0<r<r_a$. Set $r_j=d_X(x_j,a)\to0$. For all sufficiently large $j$, we have $r_j<r$. Stop each path when it first reaches background radius $r$, and parametrize this initial segment as a piecewise smooth curve $\gamma_j:[0,1]\to X_m^\circ$. Then $d_X(\gamma_j(0),a)=r_j$, $d_X(\gamma_j(1),a)=r$, and $0<d_X(\gamma_j(t),a)\le r<r_a$ throughout the segment.

By \ref{item:inverse}, for almost every $t$
$$
|\dot\gamma_j(t)|_{g_m} \ge \frac{\sqrt{c_a}}{d_X(\gamma_j(t),a)} |\dot\gamma_j(t)|_{g_X}\ge \frac{\sqrt{c_a}}{d_X(\gamma_j(t),a)}\left|\frac{d}{dt}d_X(\gamma_j(t),a)\right|.
$$
Consequently,
$$
\begin{aligned}
\mathrm{length}_{g_m}(\gamma_j) &\ge\int_0^1 \frac{\sqrt{c_a}}{d_X(\gamma_j(t),a)}\left|\frac{d}{dt}d_X(\gamma_j(t),a)\right|\,dt\\
&\ge\left|\int_0^1 \frac{\sqrt{c_a}}{d_X(\gamma_j(t),a)}\frac{d}{dt}d_X(\gamma_j(t),a)\,dt\right|\\
&=\sqrt{c_a}\log\frac{r}{r_j}\to\infty,
\end{aligned}
$$
contrary to the length bound. Hence $p_j\to a$ and continuity gives $t_m(a)=c_m$. The dimension bound follows from \cref{lem:goodlevel}.
\end{proof}
\subsection{Construction of the band}\label{subsec:construction of the band}
Let $R_m$ be the union of all connected components of $X_m^\circ$ meeting $P_m$.
\begin{proposition}\label{prop:smooth-band}
Fix $L_m>0$ sufficiently large, then there exist $W_m>0$ and a smooth function
$$
 s_m:R_m\longrightarrow\mathbb R
$$
and a closed subset
$$
 \Omega_m=s_m^{-1}([-L_m,L_m])\subset R_m
$$
with the following properties.
\begin{enumerate}[label=\textup{(B\arabic*)}]
\item\label{item:band-center}
$0,\pm L_m$ are regular values of $s_m$,
$$
 s_m^{-1}(0)=P_m,\ \sgn s_m=\sgn (t_m-c_m).
$$
\item\label{item:band-lip}
$$
 |ds_m|_{g_m}<1\text{ on }R_m.
$$

\item\label{item:band-faces}
$\Omega_m$ is a smooth manifold with boundary, with exactly the two
level-set faces
$$
 \partial\Omega_m=\partial_-\Omega_m\sqcup\partial_+\Omega_m,\ 
 \partial_\pm\Omega_m=s_m^{-1}(\pm L_m).
$$
Both faces are nonempty, and
$$
 \dist_{g_m}(\partial_-\Omega_m,\partial_+\Omega_m)\ge 2L_m.
$$

\item\label{item:band-local}
$$
 \Omega_m\subset N_{W_m}(P_m).
$$
Consequently,
$$
 \overline{\Omega_m}^{X}\cap D_m
 \subset D_m\cap t_m^{-1}(c_m),\ \dim_\HH\left(\overline{\Omega_m}^{X}\cap D_m\right)\le m-8.
$$
\end{enumerate}
\end{proposition}

\begin{proof}
Define the signed distance to $P_m$ on $R_m$ by
$$
 \rho_m=\sgn (t_m(\cdot)-c_m)\dist_{g_m}(\cdot,P_m).
$$
Then
$$
 \rho_m^{-1}(0)=P_m,\ \sgn \rho_m=\sgn (t_m-c_m),\ |\rho_m|=\dist_{g_m}(\cdot,P_m).
$$

Choose tubular neighborhood $\mathcal U_m$ and an
inner open neighborhood $\mathcal U'_m$ such that
$$
 P_m\subset\mathcal U'_m
 \subset\overline{\mathcal U'_m}^{R_m}
 \subset\mathcal U_m\subset R_m,
$$
and
$$
 \rho_m|_{\mathcal U_m}\in C^\infty(\mathcal U_m),\ 
 |d\rho_m|_{g_m}=1 \text{ on }\mathcal U_m.
$$
Choose constants
$$
 0<\eta<\frac14,\ \varepsilon_0>0,\ 
 W_m>\frac{L_m+\varepsilon_0}{1-2\eta}.
$$
and a smooth cutoff
$$
 \chi\in C^\infty(R_m;[0,1]),\ 
 \chi|_{\mathcal U'_m}=1,\ 
 \spt_{R_m}\chi\subset\mathcal U_m.
$$
Moreover
$$
 \spt_{R_m}(1-\chi)\subset R_m\setminus\mathcal U'_m,
$$
so this closed set is disjoint from $P_m$.

Choose a continuous positive function $\varepsilon:R_m\to(0,\infty)$
such that
$$
 \varepsilon(x)<\varepsilon_0,\ 
 \varepsilon(x)|d\chi(x)|_{g_m}<\frac{\eta}{2},\ \forall x\in R_m,
$$
and
$$
 \varepsilon(x)<\frac{1-2\eta}{2}|\rho_m(x)|,\ \forall x\in\spt_{R_m}(1-\chi).
$$
Apply \cite[Theorem 1]{AFLMR07}, we obtain
$$
 \widetilde h_m\in C^\infty(R_m),\ 
 |\widetilde h_m-(1-2\eta)\rho_m|<\varepsilon,\ 
 |d\widetilde h_m|_{g_m}\le1-\frac{3\eta}{2}.
$$
Now patch this approximation back to the exact signed distance near
$P_m$:
$$
 h_m=\chi(1-2\eta)\rho_m+(1-\chi)\widetilde h_m.
$$
Then $h_m$ is smooth on $R_m$, also
$$
 h_m=(1-2\eta)\rho_m \text{ on }\mathcal U'_m,\ 
 h_m=\widetilde h_m \text{ on }R_m\setminus\mathcal U_m.
$$
On $\mathcal U_m$,
$$
 dh_m=\chi(1-2\eta)d\rho_m+(1-\chi)d\widetilde h_m
       +\left((1-2\eta)\rho_m-\widetilde h_m\right)d\chi.
$$
Using the choices above,
$$
 \begin{aligned}
 |dh_m|_{g_m}
 &\le\chi(1-2\eta)+(1-\chi)(1-\frac{3\eta}{2})+\varepsilon|d\chi|_{g_m}\\
 &\le1-\eta-\chi\frac{\eta}{2}\le1-\eta.
 \end{aligned}
$$
Outside $\mathcal U_m$, the estimate follows from
$h_m=\widetilde h_m$.  Hence
$$
 |dh_m|_{g_m}\le1-\eta \text{ on }R_m.
$$
Moreover,
$$
 |h_m-(1-2\eta)\rho_m|\le(1-\chi)\varepsilon<\varepsilon_0.
$$
Therefore
$$
 h_m^{-1}(0)=P_m,\ 
 \sgn h_m=\sgn \rho_m
 =\sgn (t_m-c_m).
$$
Since $h_m=(1-2\eta)\rho_m$ near $P_m$, zero is a regular value of
$h_m$.

Choose, by Sard's theorem, a scalar $D$ such that
$$
 L_m<D<(1-2\eta)W_m-\varepsilon_0
$$
and both $D$ and $-D$ are regular values of $h_m$.  Define
$$
 s_m=\frac{L_m}{D}h_m,\ 
 \Omega_m=\left\{x\in R_m:|h_m(x)|\le D\right\}=s_m^{-1}([-L_m,L_m]).
$$
Then $0,\pm L_m$ are regular values of $s_m$, and
$$
 s_m^{-1}(0)=P_m,\ 
 \sgn s_m=\sgn (t_m-c_m),\ 
 |ds_m|_{g_m}\le\frac{L_m}{D}(1-\eta)<1.
$$
This proves \ref{item:band-center} and \ref{item:band-lip}.  

Since $\pm L_m$ are regular values, $\Omega_m$ is a smooth manifold with
boundary
$$
 \partial\Omega_m=\partial_-\Omega_m\sqcup\partial_+\Omega_m,\ 
 \partial_\pm\Omega_m=s_m^{-1}(\pm L_m).
$$
By \cref{prop:localization}, some component of $R_m$ contains points on opposite sides of $P_m$ at distances greater than $W_m$ from it. At these points the estimate $|h_m-(1-2\eta)\rho_m|<\varepsilon_0$ implies $h_m<-D$ and $h_m>D$, respectively. Hence both level sets are nonempty. For any piecewise smooth path $\gamma:[0,1]\to R_m$ from
$\partial_-\Omega_m$ to $\partial_+\Omega_m$,
$$
 2L_m
 =|s_m(\gamma(1))-s_m(\gamma(0))|
 \le\int_0^1|ds_m(\gamma(t))|_{g_m}|\dot\gamma(t)|_{g_m}\,dt
 \le\mathrm{length}_{g_m}(\gamma),
$$
and \ref{item:band-faces} follows.

Finally, if $x\in\Omega_m$, then
$$
 (1-2\eta)|\rho_m(x)|
 \le |h_m(x)|+\varepsilon_0
 \le D+\varepsilon_0
 <(1-2\eta)W_m.
$$
Hence $\dist_{g_m}(x,P_m)=|\rho_m(x)|<W_m$, which implies $\Omega_m\subset N_{W_m}(P_m)$ and \ref{item:band-local} is done.
\end{proof}
\subsection{Construction of the bubble}\label{subsec:construction of the bubble}
Write $ E_P=\left\{x\in R_m: s_m(x)<0\right\}=\left\{x\in R_m: t_m(x)<c_m\right\}$, define
$$
\mu(s)=-\frac{2\pi}{L_m}\tan\frac{\pi s}{2L_m},\ h(x)=\mu\circ s_m(x),\ \forall s\in (-L_m, L_m),\ x\in \Omega_m^\circ.
$$
Here a Caccioppoli set means a set of locally finite perimeter. For such a set $F\subset R_m$ and an open set
$V\subset\subset R_m$, write
$$
 P_{f_m}(F;V)=\int_{\partial^*F\cap V}e^{-f_m}\,d\HH^{m-1}_{g_m},
$$
where $\partial^*F$ is its reduced boundary. If $V\subset\subset\Omega_m^\circ$, define the local energy
$$
 \mathcal{F}(F;V)=P_{f_m}(F;V)-\int_{F\cap V}h e^{-f_m}\,d\HH^m_{g_m}.
$$
We will prove the following in this subsection.
\begin{proposition}[Relative weighted bubble]\label{prop:relative-weighted-bubble}
There is a Caccioppoli
set $E\subset R_m$ such that, for every open
$V\subset\subset\Omega_m^\circ$ and every Caccioppoli set
$F\subset R_m$ with $F\Delta E\subset\subset V$,
$$
 \mathcal{F}(E;V)\le\mathcal{F}(F;V).
$$
Moreover, $E$ has the following properties.
\begin{enumerate}[label=\textup{(A\arabic*)}]
\item\label{item:bv} There is an open set $U_m^\partial\subset R_m$ with
$$
 \begin{gathered}
 \partial\Omega_m\subset U_m^\partial
       \subset\left\{x\in R_m: |s_m(x)|>\frac{L_m}{2}\right\},\\
 E\setminus\Omega_m=E_P\setminus\Omega_m,\ E\cap U_m^\partial=E_P\cap U_m^\partial.
 \end{gathered}
$$
\item\label{item:support}
Its perimeter support has a regular-singular decomposition
$$
 \Sigma_m=\spt_{R_m}\HH^{m-1}_{g_m}\llcorner\partial^*E=\partial^* E\sqcup S_m\subset\Omega_m^\circ.
$$
Moreover, both $\Sigma_m$ and $S_m$ are closed in $X_m^\circ$, and $\dim_\HH S_m\le m-8$. 
\item \label{item:meancurvature}
Write $\Sigma_{m,\reg}=\partial^* E$, equip it with the induced metric $g_{\Sigma_m}=g_m|_{T\Sigma_{m,\reg}}$, and let $\nu$ be the outward unit normal of $E$ along $\Sigma_{m,\reg}$.
Then for $\varphi\in C_c^\infty(\Sigma_{m,\reg})$, one has
$$
\int_{\Sigma_{m,\reg}}\left(H-df_m(\nu)-h\right)\varphi e^{-f_m}\,d\HH^{m-1}_{g_m}=0
$$
and
$$
\int_{\Sigma_{m,\reg}}|d\varphi|_{g_{\Sigma_m}}^2e^{-f_m}\,d\HH^{m-1}_{g_m}\ge\int_{\Sigma_{m,\reg}}\left(\mathrm{Ric}_{g_m}(\nu,\nu)+|A_m|^2+D^2_{g_m}f_m(\nu,\nu)+dh(\nu)\right)\varphi^2e^{-f_m}\,d\HH^{m-1}_{g_m},
$$
where $A_m$ and $H$ are the second fundamental form and mean curvature of $\Sigma_{m,\reg}$ in $X_m^\circ$ with respect to $\nu$.
\end{enumerate}
\end{proposition}
\begin{center}
\includegraphics[width=\textwidth]{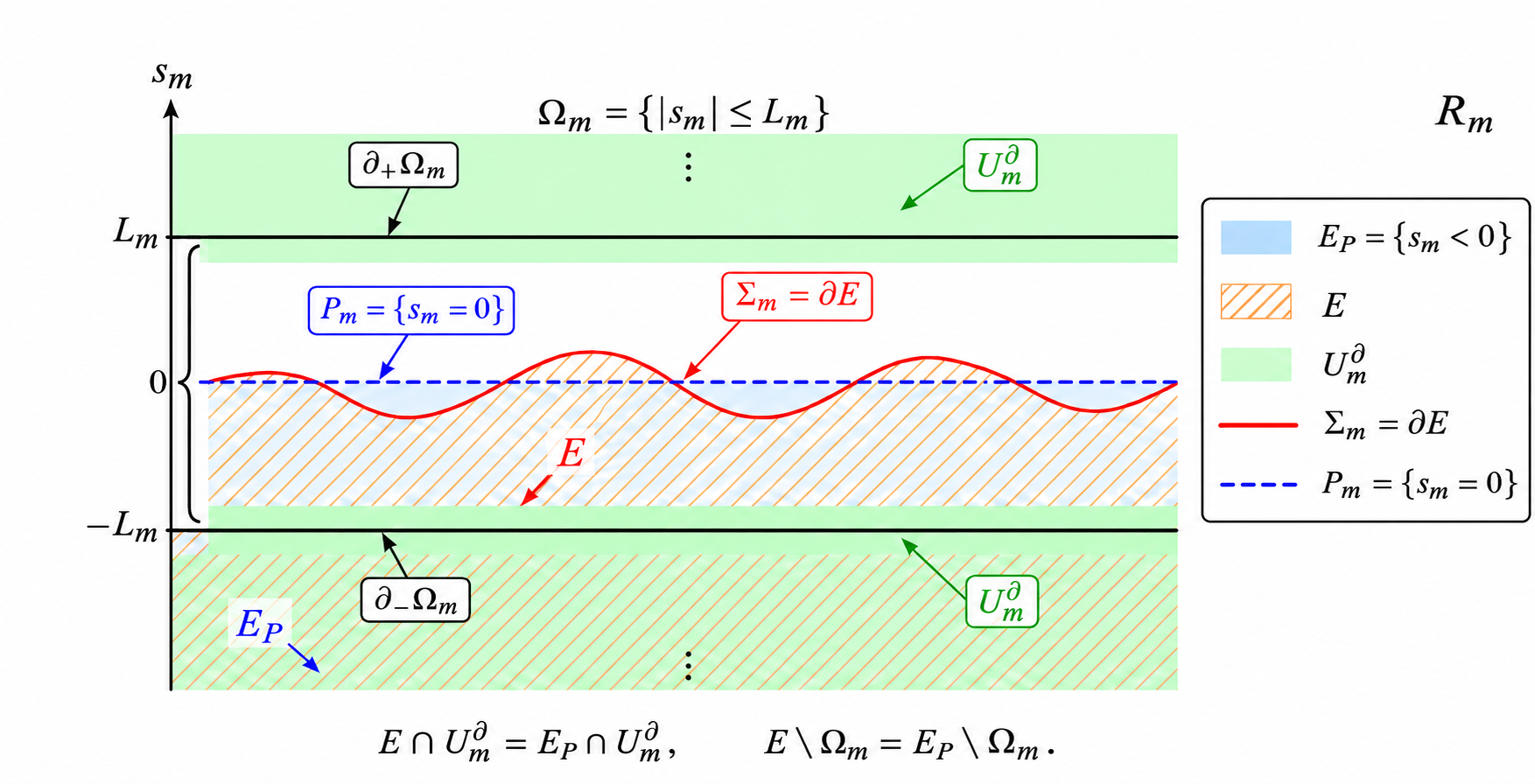}
\end{center}
\begin{proposition}\label{prop:bubble-face-barriers}
There are disjoint open Caccioppoli sets $B_-,B_+\subset R_m$, with smooth boundaries closed in $R_m$, such that
$$
 \begin{gathered}
 \left\{x\in R_m:s_m(x)\le-L_m\right\}\subset B_-\subset\left\{x\in R_m:s_m(x)<-\frac{L_m}{2}\right\},\\
 \left\{x\in R_m:s_m(x)\ge L_m\right\}\subset B_+\subset\left\{x\in R_m:s_m(x)>\frac{L_m}{2}\right\}.
 \end{gathered}
$$
The collars $B_\pm\cap\Omega_m^\circ$ carry smooth unit fields $Z_\pm$, extending across $\partial B_\pm$ and equal there to their outward normals, with
$$
 \mathrm{div}_{f_m}Z_-<h\text{ in }B_-\cap\Omega_m^\circ,\ \mathrm{div}_{f_m}Z_+<-h\text{ in }B_+\cap\Omega_m^\circ,
$$
where $\mathrm{div}_{f_m}Z_\pm=\mathrm{div}_{g_m}Z_\pm-df_m(Z_\pm)$.
\end{proposition}
\begin{proof}
Let $\nu_-$ be the inward unit normal of $\Omega_m$ along $\partial_-\Omega_m$. Choose a smooth $a: \partial_-\Omega_m\to (0, \frac{L_m}{2})$ such that $\Psi(p,t)=\exp^{g_m}_p(t\nu_-(p))$ maps $U_-=\left\{(p, t)\in \partial_-\Omega_m\times\R: 0\le t<2a(p)\right\}$ diffeomorphically onto an open neighborhood $V_-$ of $\partial_-\Omega_m$ in $\Omega_m$. 

Let $r: V_-\to \R$ be coordinate function such that $r(\Psi(p, t))=t,\ \forall (p, t)\in U_-$. For $p\in\partial_-\Omega_m$, write
$$
 M(p)=1+\max_{0\le t\le a(p)}\left.\left((m-1)\left(|D^2_{g_m}r|_{g_m}+1\right)+|df_m|_{g_m}\right)\right|_{\Psi(p,t)}.
$$
Choose a countable open cover $\left\{U_i\right\}_{i=1}^\infty$ of $\partial_-\Omega_m$, together with a subordinate locally finite smooth partition of unity $\left\{\phi_i\right\}_{i=1}^\infty$ satisfying
$$
 \phi_i\in C_c^\infty(\partial_-\Omega_m),\ \phi_i\ge0,\ \spt\phi_i\subset U_i,\ \forall i\ge 1,\ \text{and }\sum_{i=1}^\infty\phi_i=1.
$$
Extend each $\phi_i$ to $V_-$ by $\phi_i(\Psi(p,t))=\phi_i(p),\ \forall (p, t)\in U_-$. For each $i\ge 1$, choose constants 
$$
 0<c_i<\min\left\{\inf_{\spt\phi_i}\frac a2, \inf_{\spt\phi_i}\frac1M, \frac{2^{-i}}{1+\sup\limits_{\left\{\Psi(p,t):p\in\partial_-\Omega_m,\ 0\le t\le a(p)\right\}} |D^2_{g_m}\phi_i|_{g_m}}\right\}
 $$
 and define
 $$ \varepsilon=\sum_{i\ge1}c_i\phi_i.$$
Then
$$
 0<\varepsilon<\min\left\{\frac{a}2,\frac1{M}\right\},\ \sup_{\left\{\Psi(p,t):p\in\partial_-\Omega_m,\ 0\le t\le a(p)\right\}} |D^2_{g_m}\varepsilon|_{g_m}\le1.
$$
Define
$$
 B_-=\left\{x\in R_m:s_m(x)\le-L_m\right\}\cup\left\{\Psi(p, t):p\in\partial_-\Omega_m, 0<t<\varepsilon(p)\right\},\ 
 Z_-=\frac{\nabla_{g_m}(r-\varepsilon)}{|d(r-\varepsilon)|_{g_m}}\text{ on }V_-.
$$

By the Gauss lemma, $|d(r-\varepsilon)|_{g_m}^2=1+|d\varepsilon|_{g_m}^2\ge1$ on $V_-$. Therefore, on $B_-\cap \Omega_m^\circ$,
$$
 \begin{aligned}
 \mathrm{div}_{f_m}Z_-
 &=\frac{\Delta_{g_m}(r-\varepsilon)-D^2_{g_m}(r-\varepsilon)(Z_-,Z_-)}{|d(r-\varepsilon)|_{g_m}}-df_m(Z_-)\\
 &\le(m-1)\left(|D^2_{g_m}r|_{g_m}+|D^2_{g_m}\varepsilon|_{g_m}\right)+|df_m|_{g_m}\\
 &\le \frac1r<\frac{2}{s_m+L_m}\le h.
 \end{aligned}
$$
Construct similar one near $\partial_+\Omega_m$ and we are done.
\end{proof}
\begin{definition}[Face constraints]
We call a $\HH^m_{g_m}$-measurable set $E\subset R_m$ satisfies the \emph{two face constraints}
if
$$
 \HH^m_{g_m}(B_-\setminus E)=0,\ 
 \HH^m_{g_m}(E\cap B_+)=0.
$$
\end{definition}
\begin{proposition}\label{prop:bubble-unconstrained}
Let $V\subset\subset\Omega_m^\circ$ be open and let $F\subset R_m$ be a Caccioppoli set. Assume that $(B_-\setminus F)\cup(F\cap B_+)$ is contained, up to a $\HH^m_{g_m}$-null set, in a compact subset of $V$. Then the set
$$
 \widehat F=(F\cup B_-)\setminus B_+
$$
is a Caccioppoli set satisfying the two face constraints and obeys
$$
 \mathcal{F}(\widehat F;V)\le\mathcal{F}(F;V).
$$
\end{proposition}
\begin{figure}[htbp]
\centering
\includegraphics[width=\textwidth]{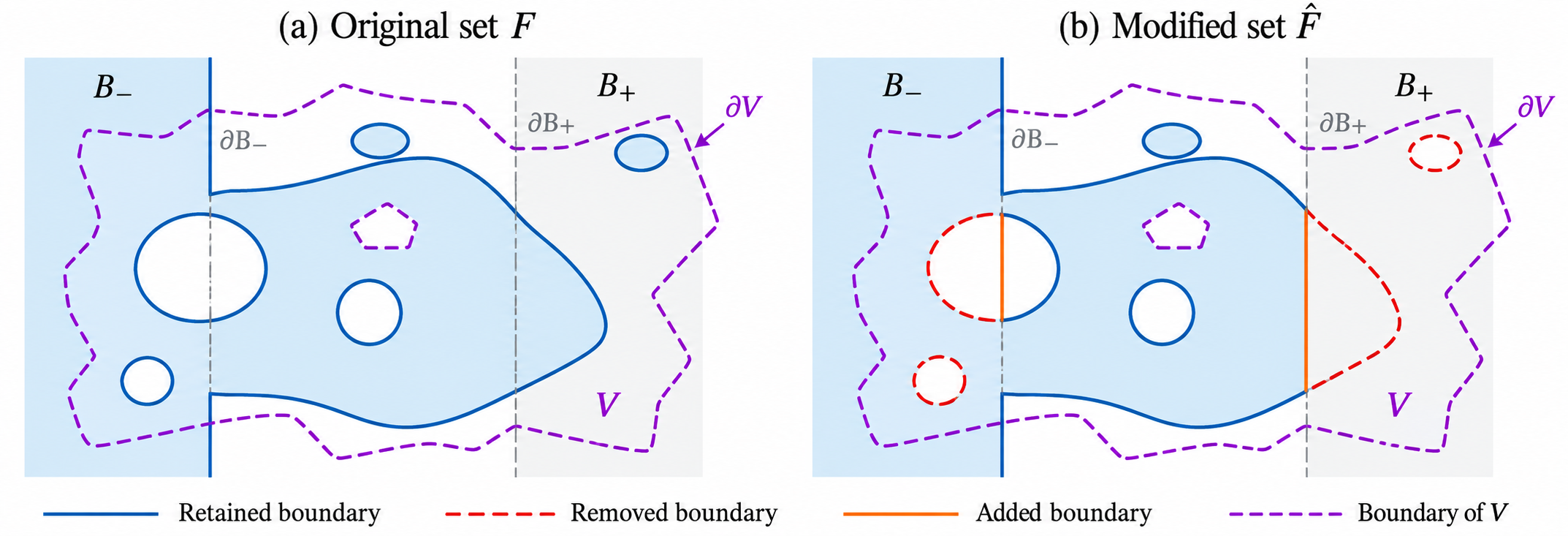}
\end{figure}
\begin{proof}
Write $F^{(0)}$ and $F^{(1)}$ for the density-zero and density-one points of $F$. By definition,
$$
 \mathcal F(\widehat F;V)-\mathcal F(F;V)
 =P_{f_m}(\widehat F;V)-P_{f_m}(F;V)-\int_{V\cap(B_-\setminus F)}h e^{-f_m}\,d\HH^m_{g_m}
 +\int_{V\cap F\cap B_+}h e^{-f_m}\,d\HH^m_{g_m}.
$$
Cancelling the unchanged boundary portions, we obtain
$$
\begin{aligned}
P_{f_m}(\widehat F;V)-P_{f_m}(F;V)
\le&\int_{V\cap\partial B_-\cap F^{(0)}}e^{-f_m}\,d\HH^{m-1}_{g_m}
+\int_{V\cap\partial B_+\cap F^{(1)}}e^{-f_m}\,d\HH^{m-1}_{g_m}\\
&-\int_{V\cap B_-\cap\partial^*F}e^{-f_m}\,d\HH^{m-1}_{g_m}
-\int_{V\cap B_+\cap\partial^*F}e^{-f_m}\,d\HH^{m-1}_{g_m}.
\end{aligned}
$$
By the Gauss--Green formula,
$$
\begin{aligned}
&P_{f_m}(\widehat F;V)-P_{f_m}(F;V)\\
&\le\int_{V\cap\partial B_-\cap F^{(0)}}\langle Z_-,\nu_{B_-}\rangle e^{-f_m}\,d\HH^{m-1}_{g_m}+\int_{V\cap\partial B_+\cap F^{(1)}}\langle Z_+,\nu_{B_+}\rangle e^{-f_m}\,d\HH^{m-1}_{g_m}\\
&-\int_{V\cap B_-\cap\partial^*F}\langle Z_-,\nu_F\rangle e^{-f_m}\,d\HH^{m-1}_{g_m}+\int_{V\cap B_+\cap\partial^*F}\langle Z_+,\nu_F\rangle e^{-f_m}\,d\HH^{m-1}_{g_m}\\
&\le\int_{V\cap\partial^*(B_-\setminus F)}\langle Z_-,\nu_{B_-\setminus F}\rangle e^{-f_m}\,d\HH^{m-1}_{g_m}+\int_{V\cap\partial^*(F\cap B_+)}\langle Z_+,\nu_{F\cap B_+}\rangle e^{-f_m}\,d\HH^{m-1}_{g_m}\\
&=\int_{V\cap(B_-\setminus F)}\mathrm{div}_{f_m}Z_-e^{-f_m}\,d\HH^m_{g_m}+\int_{V\cap F\cap B_+}\mathrm{div}_{f_m}Z_+e^{-f_m}\,d\HH^m_{g_m}.
\end{aligned}
$$
Hence, by \cref{prop:bubble-face-barriers}, 
$$
 \mathcal F(\widehat F;V)-\mathcal F(F;V)
 \le\int_{V\cap(B_-\setminus F)}(\mathrm{div}_{f_m}Z_--h)e^{-f_m}\,d\HH^m_{g_m}+\int_{V\cap F\cap B_+}(\mathrm{div}_{f_m}Z_++h)e^{-f_m}\,d\HH^m_{g_m}\le0.
$$
\end{proof}
\begin{proposition}
\label{prop:bubble-compact-stage}\label{prop:bubble-compactness}
Let $\left\{K_j\right\}_{j=1}^\infty$ be a compact exhaustion of $R_m$ with $K_j\subset\Int K_{j+1}$. Then for each $j\ge1$, there is a Caccioppoli set $E_j\subset R_m$ with $\HH^m_{g_m}\left((E_j\Delta E_P)\setminus K_j\right)=0$ satisfying the two face constraints, and
$$
 \mathcal F(E_j;V)\le\mathcal F(F;V),\ \forall V\subset\subset\Int K_j\cap\Omega_m^\circ\ \text{open},\ F\subset R_m\ \text{Caccioppoli with }F\Delta E_j\subset\subset V.
$$
Moreover, after passing to a subsequence,
$$
 \chi_{E_j}\longrightarrow\chi_E\ \hbox{in }L^1_{\mathrm{loc}}\left(R_m,\HH^m_{g_m}\right)
$$
for a Caccioppoli set $E$ satisfying the two face constraints.
\end{proposition}
\begin{proof}
Fix $j\ge1$, choose a smooth open set $V_j$ with
$$
 K_j\setminus(B_-\cup B_+)\subset V_j\subset\subset\Omega_m^\circ.
$$
Let $\mathcal A_j$ be the class of Caccioppoli sets $F\subset R_m$ satisfying the two face constraints and 
$$\HH^m_{g_m}\left((F\Delta E_P)\setminus K_j\right)=0,$$
write $I_j=\inf\limits_{F\in\mathcal A_j}\mathcal F(F;V_j).$

Since $E_P\in\mathcal A_j$ and
$$
 -\int_{V_j}|h|e^{-f_m}\,d\HH^m_{g_m}\le I_j\le\mathcal F(E_P;V_j)<\infty,
$$
for each $k\ge 1$, we can choose $F_k\in\mathcal A_j$, such that
$$
 I_j\le\mathcal F(F_k;V_j)\le I_j+\frac1k.
$$
Consequently,
$$
 \begin{aligned}
 \HH^{m-1}_{g_m}\left(\partial^*F_k\cap V_j\right)\le&\max_{\overline {V_j}}e^{f_m} P_{f_m}(F_k;V_j)\\
 \le& \max_{\overline {V_j}}e^{f_m}\left(\mathcal F(F_k;V_j)+\int_{F_k\cap V_j}h e^{-f_m}\,d\HH^m_{g_m}\right)\\
 \le&\max_{\overline {V_j}}e^{f_m}\left(\mathcal F(E_P;V_j)+1+\int_{V_j}|h|e^{-f_m}\,d\HH^m_{g_m}\right).
 \end{aligned}
$$
 Applying \cite[Theorem 1.19]{Giu84} we obtain, after passing to a subsequence, a Caccioppoli set $\hat E_j\subset V_j$ with
$$
 \chi_{F_k}\longrightarrow\chi_{\hat E_j}\ \text{in }L^1\left(V_j,\HH^m_{g_m}\right).
$$
Define $E_j=\hat E_j\cup(E_P\setminus V_j)\in\mathcal A_j$, following the argument of \cite[Theorem 1.9]{Giu84},
$$
\begin{aligned}
P_{f_m}(E_j;V_j)
&=\sup_{\substack{X\in C_c^1(V_j,TV_j)\\ |X|_{g_m}\le1}}\int_{E_j\cap V_j}\mathrm{div}_{g_m}\left(e^{-f_m}X\right)\,d\HH^m_{g_m}\\
&=\sup_{\substack{X\in C_c^1(V_j,TV_j)\\ |X|_{g_m}\le1}}\lim_{k\to\infty}\int_{F_k\cap V_j}\mathrm{div}_{g_m}\left(e^{-f_m}X\right)\,d\HH^m_{g_m}\\
&\le\varliminf_{k\to\infty}\sup_{\substack{X\in C_c^1(V_j,TV_j)\\ |X|_{g_m}\le1}}\int_{F_k\cap V_j}\mathrm{div}_{g_m}\left(e^{-f_m}X\right)\,d\HH^m_{g_m}\\
&=\varliminf_{k\to\infty}P_{f_m}(F_k;V_j).
\end{aligned}
$$
Consequently,
$$
 I_j\le\mathcal F(E_j;V_j)\le\varliminf_{k\to\infty}\mathcal F(F_k;V_j)=I_j.
$$
Let $V\subset\subset\Int K_j\cap\Omega_m^\circ$ be open and let $F\subset R_m$ be Caccioppoli with $F\Delta E_j\subset\subset V$. Write $\widehat F=(F\cup B_-)\setminus B_+$, then $\widehat F\in\mathcal A_j$ and $\widehat F\Delta E_j$ is compactly supported in $V\cap V_j$ up to null sets. By \cref{prop:bubble-unconstrained},
$$
0\le\mathcal F\left(\widehat F;V_j\right)-\mathcal F(E_j;V_j)=\mathcal F\left(\widehat F;V\right)-\mathcal F(E_j;V)\le\mathcal F(F;V)-\mathcal F(E_j;V).
$$
Now fix a smooth precompact open set $W\subset R_m$. For sufficiently large $j$, $\overline W\subset\Int K_j$ and $(E_P\cap W)\cup(E_j\setminus W)\in\mathcal A_j$, then we obtain
$$
\begin{aligned}
0&\le\mathcal F\left((E_P\cap W)\cup(E_j\setminus W);V_j\right)-\mathcal F(E_j;V_j)\\
&\le P_{f_m}(E_P;W)-P_{f_m}(E_j;W)+\int_{\partial W}e^{-f_m}\,d\HH^{m-1}_{g_m}+\int_{\overline W\setminus(B_-\cup B_+)}|h|e^{-f_m}\,d\HH^m_{g_m}.
\end{aligned}
$$
The same $BV$ compactness argument as above and a diagonal subsequence yield $\chi_{E_j}\to\chi_E$ in $L^1_{\mathrm{loc}}\left(R_m,\HH^m_{g_m}\right)$ for some Caccioppoli set $E\subset R_m$ satisfying the two face constraints. 
\end{proof}

\begin{proposition}
\label{prop:bubble-limit-minimality}\label{prop:bubble-exact-phase}
For every open $V\subset\subset\Omega_m^\circ$ and every Caccioppoli set $F$ with $F\Delta E\subset\subset V$, we have
$$
 \mathcal{F}(E;V)\le\mathcal{F}(F;V).
$$
\end{proposition}
\begin{proof}
Choose a precompact open set $W\subset R_m$ with smooth boundary such that $F\Delta E\subset\subset W\subset\subset V$ and $|D\chi_F|(\partial W)=0$. For sufficiently small $\delta>0$ and for each $0<t\le\delta$, write $W_t=\left\{x: \dist_{g_m}(x,W)<t\right\}\subset\subset V$. Coarea formula yields              
$$
\lim_{j\to \infty}\int_0^\delta\int_{\partial W_t}|\chi_{E_j}-\chi_E|e^{-f_m}\,d\HH^{m-1}_{g_m}\,dt
=\lim_{j\to \infty}\int_{W_\delta\setminus W}|\chi_{E_j}-\chi_E|e^{-f_m}\,d\HH^m_{g_m}=0.
$$
Choose $\left\{t_j\right\}_{j=1}^\infty\subset(0,\delta)$ such that
$$
\lim_{j\to \infty}\int_{\partial W_{t_j}}|\chi_{E_j}-\chi_E|e^{-f_m}\,d\HH^{m-1}_{g_m}=0.
$$
Now local minimality of $E_j$ applied to $(F\cap W_{t_j})\cup(E_j\setminus W_{t_j})$ implies
$$
\begin{aligned}
P_{f_m}(E_j;W)
&\le P_{f_m}(F;W_\delta)+\int_{\partial W_{t_j}}|\chi_{E_j}-\chi_E|e^{-f_m}\,d\HH^{m-1}_{g_m}+\int_{W_{t_j}}h(\chi_{E_j}-\chi_F)e^{-f_m}\,d\HH^m_{g_m}\\
&=P_{f_m}(F;W_\delta)+\int_{\partial W_{t_j}}|\chi_{E_j}-\chi_E|e^{-f_m}\,d\HH^{m-1}_{g_m}\\
&+\int_W h(\chi_E-\chi_F)e^{-f_m}\,d\HH^m_{g_m}+\int_{W_{t_j}}h(\chi_{E_j}-\chi_E)e^{-f_m}\,d\HH^m_{g_m}.
\end{aligned}
$$
Note that
$$
\varlimsup_{j\to\infty}\left|\int_{W_{t_j}}h(\chi_{E_j}-\chi_E)e^{-f_m}\,d\HH^m_{g_m}\right|
\le\sup_{W_\delta}|he^{-f_m}| \varlimsup_{j\to\infty}\|\chi_{E_j}-\chi_E\|_{L^1(W_\delta,\HH^m_{g_m})}=0,
$$
we obtain
$$
\begin{aligned}
P_{f_m}(E;W)
&\le\varliminf_{j\to\infty}P_{f_m}(E_j;W)\\
&\le P_{f_m}(F;W_\delta)+\int_W h(\chi_E-\chi_F)e^{-f_m}\,d\HH^m_{g_m}.
\end{aligned}
$$
Since 
$$
0\le\int_{\partial^*F\cap\partial W}e^{-f_m}\,d\HH^{m-1}_{g_m}
=\int_{\partial W}e^{-f_m}\,d|D\chi_F|\le\max_{\partial W}e^{-f_m}|D\chi_F|(\partial W)=0.
$$
By dominated convergence,
$$
\begin{aligned}
\lim_{\delta\to0}P_{f_m}(F;W_\delta)
&=\int_{\partial^*F}\lim_{\delta\to0}\chi_{W_\delta}e^{-f_m}\,d\HH^{m-1}_{g_m}\\
&=\int_{\partial^*F\cap\overline W}e^{-f_m}\,d\HH^{m-1}_{g_m}\\
&=P_{f_m}(F;W)+\int_{\partial^*F\cap\partial W}e^{-f_m}\,d\HH^{m-1}_{g_m}
=P_{f_m}(F;W).
\end{aligned}
$$
Therefore,
$$
\begin{aligned}
\mathcal F(E;V)-\mathcal F(F;V)
&=\mathcal F(E;W)-\mathcal F(F;W)\\
&=P_{f_m}(E;W)-P_{f_m}(F;W)-\int_W h(\chi_E-\chi_F)e^{-f_m}\,d\HH^m_{g_m}\le0.
\end{aligned}
$$
\end{proof}
\begin{proof}[Proof of \cref{prop:relative-weighted-bubble}]
By \cref{prop:bubble-face-barriers,prop:bubble-compact-stage,prop:bubble-limit-minimality}, choose a representative of $E$ such that $B_-\subset E$ and $E\cap B_+=\varnothing$, then \ref{item:bv} holds with
$$
 U_m^\partial=B_-\cup B_+.
$$
Fix an open set $V\subset\subset\Omega_m^\circ$, local minimality implies, for Caccioppoli set $F\subset R_m$ with $F\Delta E\subset\subset V$,
$$
P_{f_m}(E;V)\le P_{f_m}(F;V)+\int_V h(\chi_E-\chi_F)e^{-f_m}\,d\HH^m_{g_m}\le P_{f_m}(F;V)+\sup_V|h|\int_{E\Delta F}e^{-f_m}\,d\HH^m_{g_m}.
$$
Write $\widetilde g_m=e^{-\frac{2f_m}{m-1}}g_m$. Then
$$
P_{f_m}(F;V)=\HH^{m-1}_{\widetilde g_m}\left(\partial^*F\cap V\right),\ h e^{-f_m}\,d\HH^m_{g_m}=h e^{\frac{f_m}{m-1}}\,d\HH^m_{\widetilde g_m},
$$
and therefore
$$
\HH^{m-1}_{\widetilde g_m}\left(\partial^*E\cap V\right)\le\HH^{m-1}_{\widetilde g_m}\left(\partial^*F\cap V\right)+\sup_V\left|h e^{\frac{f_m}{m-1}}\right|\HH^m_{\widetilde g_m}\left(E\Delta F\right).
$$
By \cite[Theorem A.1 and Remark A.3]{Eic09}, $\partial^*E$ is $C^{1,\alpha}$, and
$$
\Sigma_m=\overline{\partial^*E}^{R_m},\ S_m=\Sigma_m\setminus\partial^*E,\ \dim_\HH S_m\le m-8.
$$
In addition, elliptic regularity makes $\partial^* E$ smooth, which proves \ref{item:support}.

For $\varphi\in C_c^\infty(\Sigma_{m,\reg})$, applying \cite[Lemma 13 and Lemma 14]{CL24} to variations compactly supported in the regular part proves \ref{item:meancurvature}.
\end{proof}
\section{Step of the descent}\label{sec:step of the descent}
Fix an $m$-dimensional datum $\cD_m$, with $m\ge 3$. We construct the next datum
$$
\cD_{m-1}=(X_{m-1}^\circ,\overline{X_{m-1}},D_{m-1},F_{m-1};g_{m-1},f_{m-1},\sigma_{m-1})
$$
by setting
$$
 X_{m-1}^\circ=\Sigma_{m,\reg},\ \overline{X_{m-1}}=\overline{\Sigma_{m,\reg}}^{X},\ D_{m-1}=\overline{X_{m-1}}\setminus\Sigma_{m,\reg},\ F_{m-1}=G_m|_{\overline{X_{m-1}}}.
$$
The metric $g_{m-1}$, weight $f_{m-1}$, and lower bound $\sigma_{m-1}$ will be constructed in \ref{subsec:next-metric}.

\subsection{Preservation of topology}\label{subsec: preservation of topology}
We prove in this subsection that the conditions \ref{item:strata}, \ref{item:map}, and \ref{item:degree} hold for $\cD_{m-1}$.
\begin{proposition}\label{prop:next-singular-set}
The set $X_{m-1}^\circ$ is an oriented smoothly embedded $(m-1)$-manifold, and $D_{m-1}$ is closed in $X$. Moreover,
$$
D_{m-1}=S_m\sqcup\left(\overline{X_{m-1}}\cap D_m\right),\ \overline{X_{m-1}}\cap D_m\subset D_m\cap t_m^{-1}(c_m),\ \dim_\HH D_{m-1}\le m-8.
$$
\end{proposition}
\begin{proof}
By \ref{item:support}, $\Sigma_{m,\reg}$ is smooth and relatively open in $\Sigma_m$, and $\Sigma_m$ is closed in $X_m^\circ$. Consequently,
$$
\overline{X_{m-1}}\cap X_m^\circ=\overline{\Sigma_{m,\reg}}^{X_m^\circ}=\Sigma_m=\Sigma_{m,\reg}\sqcup S_m.
$$
Since $\overline{X_{m-1}}\subset\overline{X_m}=X_m^\circ\sqcup D_m$, removing $\Sigma_{m,\reg}$ yields the decomposition. The regular part is relatively open in $\overline{X_{m-1}}$, so $D_{m-1}$ is closed in $X$. By \ref{item:band-local} and \cref{lem:goodlevel},
$$
\overline{X_{m-1}}\cap D_m\subset\overline{\Omega_m}^{X}\cap D_m\subset D_m\cap t_m^{-1}(c_m),\ \dim_\HH\left(D_m\cap t_m^{-1}(c_m)\right)\le m-8.
$$
Together with $\dim_\HH S_m\le m-8$, \ref{item:strata} is done.
\end{proof}

\begin{proposition}\label{prop:next-proper-map}
The map $G_m$ is proper on $\overline{\Omega_m}^{X}$. In particular, $F_{m-1}$ is proper, locally Lipschitz, and smooth on $X_{m-1}^\circ$.
\end{proposition}
\begin{proof}
Write $\Lambda_m=\sup\limits_X|d\Phi_{m-2}|_{g_X}$. Since $g_m\ge g_X$ and $\Omega_m\subset N_{W_m}(P_m)$,
$$
|t_m(x)-c_m|\le \Lambda_m\dist_{g_m}(x,P_m)\le\Lambda_m W_m,\ \forall x\in \Omega_m.
$$
For compact $K\subset\R^{m-3}$,
$$
\overline{\Omega_m}^{X}\cap G_m^{-1}(K)\subset F_m^{-1}\left(K\times[c_m-\Lambda_m W_m,c_m+\Lambda_m W_m]\right),
$$
and the set on the left is closed in the compact set on the right by \ref{item:map}. Since $\overline{X_{m-1}}$ is closed in $\overline{\Omega_m}^{X}$, its restriction $F_{m-1}$ is proper. Local Lipschitz continuity and smoothness follow by restriction from $G_m$, proving \ref{item:map}.
\end{proof}

Write
$V_{m-1}=\R^{m-3}\setminus F_{m-1}(D_{m-1}),\ N_{m-1}=F_{m-1}^{-1}(V_{m-1}),\ U_m=V_{m-1}\setminus G_m\left(D_m\cap t_m^{-1}(c_m)\right).$
\begin{proposition}\label{prop:detected-degree-transfer}
The map $F_{m-1}|_{N_{m-1}}:N_{m-1}\to V_{m-1}$ is surjective and
$$
\udeg(F_{m-1}|_{N_{m-1}})=\udeg(F_m|_{N_m}).
$$
More precisely, for a.e. $z\in U_m$,
the fibres
$$
\Sigma_{m,z}=\Sigma_{m,\reg}\cap G_m^{-1}(z),\ P_{m,z}=P_m\cap G_m^{-1}(z)
$$
are compact oriented smooth surfaces and
$$
[\Sigma_{m,z}]=[P_{m,z}]\ \text{in }H_2(X_m^\circ;\Z).
$$
\end{proposition}
\begin{proof}
Write $\llbracket E\rrbracket$ and $\llbracket E_P\rrbracket$ for their integration currents in the oriented manifold $X_m^\circ$, and use the same notation for the oriented hypersurfaces $\Sigma_{m,\reg}$ and $P_m$. Define
$$
\cB_m=(-1)^{m-1}\left(\llbracket E\rrbracket-\llbracket E_P\rrbracket\right).
$$
Then Gauss--Green formula, \ref{item:bv} and \ref{item:support} imply
$$
\partial\cB_m=\llbracket\Sigma_{m,\reg}\rrbracket-\llbracket P_m\rrbracket,\ \spt\cB_m\subset\Omega_m.
$$
In addition, for every open $V\subset\subset X_m^\circ$, the mass in $V$ satisfies
$$
\mathbf M_V(\cB_m)=\int_V|\chi_E-\chi_{E_P}|\,d\HH^m_{g_m}
=\HH^m_{g_m}\left((E\Delta E_P)\cap V\right)\le\HH^m_{g_m}(V)<\infty,$$
$$
\mathbf M_V(\partial\cB_m)\le\HH^{m-1}_{g_m}(\partial^*E\cap V)+\HH^{m-1}_{g_m}(P_m\cap V)<\infty.
$$
Therefore, we obtain that $\cB_m$ is locally integral.

Suppose first that $m>3$, since $\cB_m$ is locally integral, $G_m$ is locally Lipschitz, and $\spt\cB_m$ is separable, \cite[Theorems 6.4 and 8.5(1)]{Lan11} ensures that, for a.e. $z\in U_m$,
$$
\cB_{m,z}=(-1)^{m-3}\langle\cB_m,G_m,z\rangle
$$
exists as a locally integral three-current. Explicitly, for every compactly supported smooth three-form $\omega$,
$$
\cB_{m,z}(\omega)
=\lim_{\varepsilon\to 0}\frac{1}{\varepsilon^{m-3}}
\int_{X_m^\circ\cap G_m^{-1}\left(\prod\limits_{i=1}^{m-3}[z_i,z_i+\varepsilon]\right)}
(\chi_E-\chi_{E_P})\,d\Phi_1\wedge\cdots\wedge d\Phi_{m-3}\wedge\omega;
$$
see \cite[Definition 6.3 and equation (6.14)]{Lan11}.

Choose $z$ outside the null exceptional sets for these slices, for the boundary slices, for \cref{lem:goodlevel}, and for the critical values of $G_m|_{P_m}$ and $G_m|_{\Sigma_{m,\reg}}$. The two endpoint fibres $\Sigma_{m, z}$ and $P_{m, z}$ are then smooth oriented surfaces. 

By \cite[equation (6.8)]{Lan11} and the compactness of $\overline{\Omega_m}^{X}\cap G_m^{-1}(z)$, we have
$$
\spt\cB_{m,z}\subset\spt\cB_m\cap G_m^{-1}(z)\subset\overline{\Omega_m}^{X}\cap G_m^{-1}(z)\subset\subset X_m^\circ.
$$
Thus $\cB_{m,z}$ has compact support, finite mass and boundary mass. 

The slicing boundary formula \cite[equation (6.9)]{Lan11} now yields
$$
\begin{aligned}
\partial\cB_{m,z}
&=(-1)^{m-3}\partial\langle\cB_m,G_m,z\rangle\\
&=\langle\partial\cB_m,G_m,z\rangle\\
&=\llbracket\Sigma_{m,z}\rrbracket-\llbracket P_{m,z}\rrbracket.
\end{aligned}
$$
For $m=3$, $U_m=\R^0=\{\mathrm{pt}\}$ and no slicing is required. Write
$$
\cB_{3,\mathrm{pt}}=\cB_3,\ \cB_{3,\mathrm{pt}}(\omega)=\int_{X_3^\circ}(\chi_E-\chi_{E_P})\omega,
$$
and the above properties hold.

By \ref{item:metric}, each component of $(X_m^\circ,g_m)$ is complete and admits locally strong Lipschitz contractions in bi-Lipschitz charts. Since $\cB_{m,z}$ has compact support, the integral-current and singular homology comparison \cite[Corollaries 1.4 and 1.6]{RS09}, applied componentwise, yields
$$
[\Sigma_{m,z}]=[P_{m,z}]\ \text{in }H_2(X_m^\circ;\Z).
$$
Pairing with $u_X|_{X_m^\circ}$ and using \cref{lem:goodlevel}, we obtain $\langle u_X|_{\Sigma_{m,z}},[\Sigma_{m,z}]\rangle=\udeg(F_m|_{N_m})\ne0$. Thus the image of $F_{m-1}|_{N_{m-1}}$ contains a full-measure, dense subset of $V_{m-1}$ and is closed by properness, hence equals $V_{m-1}$. Now \cref{prop:fibre degree} implies
$$
\udeg(F_{m-1}|_{N_{m-1}})=\langle u_X|_{\Sigma_{m,z}},[\Sigma_{m,z}]\rangle=\udeg(F_m|_{N_m}),
$$
which is \ref{item:degree} for the next stage.
\end{proof}
\subsection{Construction of the next metric}\label{subsec:next-metric}
Following \cite[equation (2.10)]{BZ26}, let $(N^d,g)$ be a smooth Riemannian manifold, $f\in C^\infty(N)$, and $a,\sigma\in\R$. For $\varphi\in C_c^\infty(N)$, define
$$
Q_{a,\sigma}^{g,f}(\varphi)=\int_N|d\varphi|_g^2e^{-f}\,d\HH^d_g+\frac12\int_N\left(R_g+2\Delta_g f-a|df|_g^2-\sigma\right)\varphi^2e^{-f}\,d\HH^d_g.
$$
We use this form on $\Sigma_{m,\reg}$ and its connected components. In hypersurface expressions, $f_m$ denotes its restriction to $\Sigma_{m,\reg}$. Write
$$
a=\frac{m}{m+1},\ \tau_m=\sigma_m-\frac{2\pi^2}{L_m^2}.
$$
We take $L_m$ large enough in \cref{prop:smooth-band} so that $\tau_m>0$.
\begin{proposition}\label{prop:inherited-stability}
For every $\varphi\in C_c^\infty(\Sigma_{m,\reg})$,
$$
Q_{a,\tau_m}^{g_{\Sigma_m},f_m}(\varphi)\ge0.
$$
\end{proposition}
\begin{proof}
Note that
$$
\begin{aligned}
2\mathrm{Ric}_{g_m}(\nu,\nu)&=R_{g_m}-R_{g_{\Sigma_m}}+H^2-|A_m|^2,\\
\Delta_{g_m}f_m&=\Delta_{g_{\Sigma_m}}f_m+D^2_{g_m}f_m(\nu,\nu)+Hdf_m(\nu),\\
|df_m|_{g_m}^2&=|df_m|_{g_{\Sigma_m}}^2+df_m(\nu)^2.
\end{aligned}
$$
we obtain
$$
\begin{aligned}
&\mathrm{Ric}_{g_m}(\nu,\nu)+|A_m|^2+D^2_{g_m}f_m(\nu,\nu)+dh(\nu)\\
&=\frac12\left(R_{g_m}-R_{g_{\Sigma_m}}+H^2+|A_m|^2\right)
+D^2_{g_m}f_m(\nu,\nu)+dh(\nu)\\
&=\frac12\left[\mathcal S_{f_m}^{m,-1}(g_m)-R_{g_{\Sigma_m}}+a|df_m|_{g_m}^2\right]
-\Delta_{g_m}f_m+D^2_{g_m}f_m(\nu,\nu)\\
&\ +\frac12H^2+\frac12|A_m|^2+dh(\nu)\\
&=\frac12\left[\mathcal S_{f_m}^{m,-1}(g_m)-R_{g_{\Sigma_m}}-2\Delta_{g_{\Sigma_m}}f_m+a|df_m|_{g_{\Sigma_m}}^2\right]\\
&\ +\frac a2df_m(\nu)^2+\frac12|A_m|^2+\frac12H^2-Hdf_m(\nu)+dh(\nu).
\end{aligned}
$$
By Cauchy--Schwarz,
$$
H^2=(\mathrm{tr}_{g_{\Sigma_m}}A_m)^2\le(m-1)|A_m|^2.
$$
Consequently,
$$
\begin{aligned}
&\frac a2df_m(\nu)^2+\frac12|A_m|^2+\frac12H^2-Hdf_m(\nu)\\
&\ge\frac{m}{2(m+1)}df_m(\nu)^2+\frac{m}{2(m-1)}H^2-Hdf_m(\nu)\\
&=\frac{m}{2(m+1)}df_m(\nu)^2+\frac{m}{2(m-1)}\left(df_m(\nu)+h\right)^2
-\left(df_m(\nu)+h\right)df_m(\nu)\\
&=\frac1{m^2-1}\left(df_m(\nu)+\frac{m+1}{2}h\right)^2+\frac14h^2
\ge\frac14h^2.
\end{aligned}
$$
Since 
$$
\mu(s)=-\frac{2\pi}{L_m}\tan\frac{\pi s}{2L_m},\ \mu'(s)=-\frac{\pi^2}{L_m^2}\cos^{-2}\frac{\pi s}{2L_m}<0,\ ds_m(\nu)\le|ds_m|_{g_m}|\nu|_{g_m}<1.
$$
It follows that
$$
\begin{aligned}
2dh(\nu)+\frac12h^2
&=2\mu'(s_m)ds_m(\nu)+\frac12\mu(s_m)^2\\
&\ge2\mu'(s_m)+\frac12\mu(s_m)^2\\
&=-\frac{2\pi^2}{L_m^2}\cos^{-2}\frac{\pi s_m}{2L_m}
+\frac{2\pi^2}{L_m^2}\tan^2\frac{\pi s_m}{2L_m}
=-\frac{2\pi^2}{L_m^2}.
\end{aligned}
$$
Combining these estimates with $\mathcal S_{f_m}^{m,-1}(g_m)\ge\sigma_m$ and the stability inequality in \ref{item:meancurvature}, we conclude that
$$
\begin{aligned}
&\int_{\Sigma_{m,\reg}}|d\varphi|_{g_{\Sigma_m}}^2e^{-f_m}\,d\HH^{m-1}_{g_m}\\
&\ge\int_{\Sigma_{m,\reg}}\left(\mathrm{Ric}_{g_m}(\nu,\nu)+|A_m|^2
+D^2_{g_m}f_m(\nu,\nu)+dh(\nu)\right)\varphi^2e^{-f_m}\,d\HH^{m-1}_{g_m}\\
&\ge\frac12\int_{\Sigma_{m,\reg}}\left[\sigma_m-\frac{2\pi^2}{L_m^2}
-R_{g_{\Sigma_m}}-2\Delta_{g_{\Sigma_m}}f_m+a|df_m|_{g_{\Sigma_m}}^2\right]\varphi^2e^{-f_m}\,d\HH^{m-1}_{g_m}.
\end{aligned}
$$
\end{proof}

\begin{lemma}\label{lem:singular-blowup-potential}
Suppose $S_m\ne\varnothing$, so $m\ge8$ by the bound $\dim_\HH S_m\le m-8$. For $u\in C^\infty(\Sigma_{m,\reg})$, write
$$
\alpha=\frac{m+1}{m-3},\ \gamma=\frac{m-3}{m+1},\ L_\gamma u=\Delta_{g_{\Sigma_m}}u-\gamma\langle df_m,du\rangle_{g_{\Sigma_m}}.
$$
There is a nonnegative $\Theta\in C^\infty(\Sigma_{m,\reg})$ such that
$$
-L_\gamma\Theta\ge\Theta^{1+\alpha}-1.
$$
For every $a_0\in S_m$, there are $c_{a_0},r_{a_0}>0$ such that
$$
\Theta(x)\ge c_{a_0}\dist_{g_m}(x,S_m)^{-\frac{2}{\alpha}},\ \forall x\in\Sigma_{m,\reg}\cap B_{g_m}(a_0,r_{a_0}).
$$
\end{lemma}
\begin{proof}
Since $E$ is locally minimizing and $\Sigma_m\subset\Omega_m^\circ$, the interior argument of \cite[Proposition 4.6]{BZ26} establishes the local AM--PI structure and packing codimension at least seven.

Apply \cite[Proposition 3.8]{BZ26} in dimension $m-1$, with singular set $S_m$, weight $e^{-\gamma f_m}$, and growth exponent
$$
\frac2\alpha=\frac{2(m-3)}{m+1}<2<7-2.
$$
This supplies the required potential and growth estimate. Multiplying the potential by a sufficiently small positive constant normalizes its differential inequality to $-L_\gamma\Theta\ge\Theta^{1+\alpha}-1$ and preserves its growth. Take $\Theta=0$ on components of $\Sigma_m$ disjoint from $S_m$.
\end{proof}

\begin{proposition}\label{prop:next-metric}
For every $0<\delta_m<\tau_m$, there are smooth $g_{m-1}$ and $f_{m-1}$ on $X_{m-1}^\circ$ satisfying \ref{item:metric} and \ref{item:inverse}, with
$$
g_{m-1}\ge g_{\Sigma_m},\ \sigma_{m-1}=\tau_m-\delta_m>0.
$$
If $S_m=\varnothing$, one may omit $\delta_m$ and take $g_{m-1}=g_{\Sigma_m}$.
\end{proposition}
\begin{proof}
Assume $S_m\ne\varnothing$. Choose
$$
0<\varepsilon\le\frac{\delta_m}{\alpha\left(1+|\tau_m|^{1+\frac 1\alpha}\right)}
$$
and define
$$
w=1+\varepsilon\Theta,\ g_{m-1}=w^\alpha g_{\Sigma_m},\ \widetilde f_m=f_m+\frac{m+1}{2}\log w.
$$
Applying \cite[Proposition 2.3]{BZ26} in dimension $m-1$ yields
$$
Q_{a,\tau_m-\delta_m}^{g_{m-1},\widetilde f_m}(\varphi)
=Q_{a,\tau_m}^{g_{\Sigma_m},f_m}(\varphi)+\frac12\int_{\Sigma_{m,\reg}}\left[-\alpha \frac{L_\gamma w}{w}-\tau_m(w^\alpha-1)+\delta_mw^\alpha\right]\varphi^2e^{-f_m}\,d\HH^{m-1}_{g_m}.
$$
Since $w\ge1$, $-L_\gamma w\ge\varepsilon(\Theta^{1+\alpha}-1)$, 
$$
w^\alpha-1=\int_{1}^{w}\alpha t^{\alpha-1}\,dt\le\int_{1}^{w}\alpha w^{\alpha-1}\,dt=\alpha(w-1)w^{\alpha-1}=\alpha\varepsilon\Theta w^{\alpha-1},$$
and
$$|\tau_m|\frac\Theta w\le\frac {1}{\alpha+1}\left(\frac\Theta w\right)^{1+\alpha}+\frac{\alpha}{1+\alpha}|\tau_m|^{1+\frac 1\alpha}\le\left(\frac\Theta w\right)^{1+\alpha}+|\tau_m|^{1+\frac 1\alpha},
$$
we have
$$
\begin{aligned}
&-\alpha\frac{L_\gamma w}{w}-\tau_m(w^\alpha-1)+\delta_mw^\alpha\\
&\ge\frac{\alpha\varepsilon\Theta^{1+\alpha}}w-\frac{\alpha\varepsilon}w
-\alpha|\tau_m|\varepsilon\Theta w^{\alpha-1}+\delta_mw^\alpha\\
&=w^\alpha\left[\alpha\varepsilon\left(\frac\Theta w\right)^{1+\alpha}
-\frac{\alpha\varepsilon}{w^{1+\alpha}}-\alpha|\tau_m|\varepsilon\frac\Theta w+\delta_m\right]\\
&\ge w^\alpha\left[\delta_m-\alpha\varepsilon\left(1+|\tau_m|^{1+\frac 1\alpha}\right)\right]\ge0.
\end{aligned}
$$
Therefore \cref{prop:inherited-stability} implies
$$
Q_{a,\sigma_{m-1}}^{g_{m-1},\widetilde f_m}(\varphi)
\ge Q_{a,\tau_m}^{g_{\Sigma_m},f_m}(\varphi)\ge0.
$$
If $S_m=\varnothing$, let $g_{m-1}=g_{\Sigma_m}$, $\widetilde f_m=f_m$, and $\sigma_{m-1}=\tau_m$ and the same conclusion follows directly from \cref{prop:inherited-stability}.

We construct a smooth function $\psi\in C^\infty(X_{m-1}^\circ)$ satisfying $\psi>0$ and $\mathcal L\psi\ge0$, where $\mathcal L:C^\infty(X_{m-1}^\circ)\to C^\infty(X_{m-1}^\circ)$ is defined by
$$
\mathcal L\psi=-\Delta_{g_{m-1}}\psi+\langle d\widetilde f_m,d\psi\rangle_{g_{m-1}}
+\frac12\left(R_{g_{m-1}}+2\Delta_{g_{m-1}}\widetilde f_m-a|d\widetilde f_m|_{g_{m-1}}^2-\sigma_{m-1}\right)\psi.
$$
Fix a connected component $N\subset X_{m-1}^\circ$ and $x_0\in N$. Integration by parts and the preceding quadratic-form inequality imply, for $\varphi\in C_c^\infty(N)$,
$$
\int_N\varphi\mathcal L\varphi\,e^{-\widetilde f_m}\,d\HH^{m-1}_{g_{m-1}}
=Q_{a,\sigma_{m-1}}^{g_{m-1},\widetilde f_m}(\varphi)\ge0.
$$
If $N$ is compact, its first eigenvalue is nonnegative and its first eigenfunction can be chosen strictly positive, so $\mathcal L\psi\ge0$.

Suppose $N$ is noncompact. For $j\ge 1$, inductively choose connected open set $U_j\subset\subset N$ with $\overline{U_{j}}\subset U_{j+1}$ and smooth boundary such that
$$
x_0\in U_1,\ \bigcup_{j=1}^\infty U_j=N.
$$
The first Dirichlet eigenvalue on $U_j$ is
$$
\lambda_j=\inf_{0\ne\varphi\in C_c^\infty(U_j)}
\frac{Q_{a,\sigma_{m-1}}^{g_{m-1},\widetilde f_m}(\varphi)}
{\displaystyle\int_{U_j}\varphi^2e^{-\widetilde f_m}\,d\HH^{m-1}_{g_{m-1}}}\ge0.
$$
Choose a first eigenfunction $\psi_j\in C^\infty(U_j)\cap C^0(\overline{U_j})$, positive on $U_j$, such that
$$
\mathcal L\psi_j=\lambda_j\psi_j\text{ in }U_j,\ \psi_j=0\text{ on }\partial U_j,\ \psi_j(x_0)=1.
$$
Domain inclusion in the Rayleigh quotient implies
$$
0\le\lambda_{j+1}\le\lambda_j\le\lambda_1,\ \lambda_j\longrightarrow\lambda\ge0.
$$
For each compact $K\subset N$,  Harnack's inequality implies for some $C_K>0$,
$$
C_K^{-1}\le\psi_j\le C_K\text{ on }K,
$$
for all sufficiently large $j$. 

Interior elliptic estimates then imply, for every integer $k\ge0$, and some $C_{K, k}>0$
$$
\|\psi_j\|_{C^k(K)}\le C_{K,k}
$$
for all sufficiently large $j$. Hence up to a subsequence, one has
$$
\psi_j\longrightarrow\psi\text{ in }C^\infty_{\mathrm{loc}}(N),\ \psi(x_0)=1,\ \mathcal L\psi=\lambda\psi.
$$
Passing to the limit in the Harnack bounds shows $\psi>0$ on $N$, and $\lambda\ge0$ implies $\mathcal L\psi\ge0$.

The functions constructed on the connected components define a smooth positive $\psi$ on $X_{m-1}^\circ$ satisfying $\mathcal L\psi\ge0$. Write $f_{m-1}=\widetilde f_m-\log\psi$, then
$$
\mathcal S_{f_{m-1}}^{m-1,-1}(g_{m-1})
=\sigma_{m-1}+2\frac{\mathcal L\psi}{\psi}
+\frac{m+1}{m}\left|d\log\psi-\frac1{m+1}d\widetilde f_m\right|_{g_{m-1}}^2
\ge\sigma_{m-1}.
$$

By construction, $g_{m-1}\ge g_{\Sigma_m}\ge g_X$. If $a_0\in S_m$, let $c_{a_0}>0$ be the constant in \cref{lem:singular-blowup-potential}. Since $w=1+\varepsilon\Theta\ge\varepsilon\Theta$, and by smooth metric comparison in an embedded chart of $X_m^\circ$ near $a_0$, there is $C_{a_0}>0$ such that
$$
w^\alpha\ge(\varepsilon c_{a_0})^\alpha\dist_{g_m}(x,S_m)^{-2},\ \dist_{g_m}(x,S_m)\le\dist_{g_m}(x,a_0)\le C_{a_0}d_X(x,a_0).
$$
Consequently, 
$$
g_{m-1}=w^\alpha g_{\Sigma_m}
\ge\frac{(\varepsilon c_{a_0})^\alpha}{\dist_{g_m}(x,S_m)^2}g_X
\ge\frac{(\varepsilon c_{a_0})^\alpha}{C_{a_0}^2d_X(x,a_0)^2}g_X.
$$
At $a_0\in\overline{X_{m-1}}\cap D_m$, an estimate of the required form follows by restricting \ref{item:inverse} to $T X_{m-1}^\circ$. Thus \ref{item:inverse} holds at every point of $D_{m-1}$.

Finally, let $\left\{x_j\right\}_{j=1}^\infty$ be a $g_{m-1}$-Cauchy sequence in one component. By \ref{item:metric} and \ref{item:support},
$$
\varlimsup_{i, j\to\infty}d_{g_m}(x_i,x_j)\le \varlimsup_{i, j\to\infty}d_{g_{m-1}}(x_i,x_j)=0.
$$
Completeness of the ambient component and closedness of $\Sigma_m$ in $X_m^\circ$ imply that $x_j$ converges in $g_m$ to some $a_0\in\Sigma_m$. If $a_0\in S_m$, choose $0<r<d_X(x_1,a_0)$ sufficiently small. Then
$$
\varliminf_{j\to\infty}d_{g_{m-1}}(x_j,x_1)\ge\varliminf_{j\to\infty}\frac{(\varepsilon c_{a_0})^{\frac{\alpha}{2}}}{C_{a_0}}\log\frac{r}{d_X(x_j,a_0)}=\infty,
$$
contradicting the Cauchy property. Thus $a_0\in\Sigma_{m,\reg}$. Near this regular point the embedded hypersurface topology agrees with the topology induced by the smooth metric $g_{m-1}$. Hence $x_j\to a_0$ in $g_{m-1}$, proving \ref{item:metric}.
\end{proof}

\subsection{The two-dimensional obstruction}\label{subsec:completion-of-descent}
\begin{proposition}\label{prop:datum-obstruction}
There is no descent datum $\cD_m$, $m\ge2$, satisfying \cref{def:datum}.
\end{proposition}
\begin{proof}
For $m\ge3$, iterate \cref{sec:step of the descent}, choosing $L_j$ and $\delta_j$ so that 
$$\sum_{j=3}^m\bigl(\frac{2\pi^2}{L_j^2}+\delta_j\bigr)<\sigma_m
.$$ 
If $m=2$, use the given datum. In either case the terminal datum $\cD_2$ satisfies $\sigma_2>0$ and has nonzero detected degree.

For the terminal datum $\cD_2$, condition \ref{item:strata} implies $D_2=\varnothing$. Hence $F_2:X_2^\circ\to\R^0$ is proper by \ref{item:map}, so $X_2^\circ$ is a closed oriented surface. Write $X_2^\circ=C_1\sqcup\cdots\sqcup C_r$ for its connected components. On each component $C$, condition \ref{item:metric} implies
$$
\begin{aligned}
4\pi\chi(C)&=\int_C R_{g_2}\,d\HH^2_{g_2}\\
&\ge\sigma_2\HH^2_{g_2}(C)+\frac23\int_C|df_2|_{g_2}^2\,d\HH^2_{g_2}>0.
\end{aligned}
$$
 Thus every component of $X_2^\circ$ is diffeomorphic to $S^2$. Since $u_X$ annihilates spherical classes in $X$,
$$
\langle u_X|_C,[C]_\Z\rangle=0.
$$
The degree formula and condition \ref{item:degree} now imply
$$
0=\sum_{i=1}^r\langle u_X|_{C_i},[C_i]_\Z\rangle=\udeg(F_2)\ne0,
$$
a contradiction.
\end{proof}
\begin{proof}[Proof of \cref{thm:abstract-obstruction}]
If $\udeg(\Phi)\ne0$, then $(X,X,\varnothing,\Phi;g_X,0,\sigma)$ is a descent datum in dimension $n$, contrary to \cref{prop:datum-obstruction}.
\end{proof}

\appendix
\section{A width bound for SYS bands}\label{app:sys-width}

Gromov asked whether the constant in his width estimate for SYS bands could be reduced to $2\pi\sqrt{\frac{n-1}{\sigma n}}$ \cite[Section~5]{Gro18}. Here we use one separating $\mu$-bubble and the descent obstruction of \cref{prop:datum-obstruction} to obtain the estimate in all dimensions.

\begin{definition}[Proper band]\label{def:band}
A proper band is a Riemannian manifold $(V,g)$ with a specified decomposition $\partial V=\partial_-V\sqcup\partial_+V$, where each of $\partial_\pm V$ is a nonempty union of connected components of $\partial V$. The width $\width(V,g)$ is the intrinsic distance between $\partial_-V$ and $\partial_+V$.
\end{definition}
\begin{definition}[Band map]\label{def:band map}
A smooth map $F: V_1\to V_2$ between two proper bands is called a band map if 
$$F(\partial_\pm V_1)\subset \partial_\pm V_2.$$
\end{definition}
\begin{definition}[SYS band]\label{def:sys-band}
A compact oriented proper band $V^n$ is an SYS band if there is a closed oriented SYS manifold $Y^{n-1}$ and a band map $F:V\to Y\times[-1,1]$ with relative degree 1, i.e.
$$
F_*[V,\partial V]_{\Z}=[Y\times[-1,1],Y\times\qty{-1,1}]_{\Z}\ \text{in } H_n(Y\times[-1,1],Y\times\qty{-1,1};\Z).
$$
\end{definition}

\begin{theorem}[SYS band inequality]\label{prop:sys-width}
If a compact oriented SYS band $(V^n,g)$ satisfies $R_g\ge\sigma>0$, then
$$
\width(V,g)\le 2\pi\sqrt{\frac{n-1}{n\sigma}}.
$$
The constant is optimal for every fixed SYS manifold $Y$ in \cref{def:sys-band}.
\end{theorem}

\begin{lemma}[A stable separator]\label{lem:band-separator}
Let $(V^n,g)$ be a compact oriented proper band of dimension $n\ge2$ with $R_g\ge\sigma$. If $0<\ell<\width(V,g)$, then there exists a Caccioppoli set $E\subset V$ containing a collar of $\partial_-V$ and disjoint from a collar of $\partial_+V$. Its interior perimeter support $\Sigma$ is compact in $\Int V$, with smooth regular part $\Sigma_{\reg}$ and singular set $D$ of Hausdorff dimension at most $n-8$. For every $\varphi\in C_c^\infty(\Sigma_{\reg})$,
$$
\int_{\Sigma_{\reg}}\bigl(2|d\varphi|^2+R_{\Sigma_{\reg}}\varphi^2\bigr)\,d\HH^{n-1}_g
\ge\left(\sigma-\frac{4\pi^2(n-1)}{n\ell^2}\right)\int_{\Sigma_{\reg}}\varphi^2\,d\HH^{n-1}_g.
$$
\end{lemma}
\begin{proof}
Choose $\ell<2a<\width(V,g)$. By \cite[Lemma~7.2]{CZ24}, there is a smooth $1$-Lipschitz function $s:V\to[-a,a]$ with $s=\pm a$ on $\partial_\pm V$. Choose $\ell<L<2a$ so that both $\pm\frac{L}{2}$ are regular values of $s$. On the interior of $U=\{|s|\le\frac{L}{2}\}$, set
$$
h(s)=-\frac{2\pi(n-1)}{nL}\tan\frac{\pi s}{L}.
$$
 On each component of $U$ meeting both faces, apply \cite[Proposition~6.12]{Cho21} with the lower face as the positive face in that proposition. The resulting minimizing Caccioppoli set contains a collar of the lower face and is disjoint from a collar of the upper face. On any component of $U$ meeting only the lower face take the whole component, and on one meeting only the upper face or neither face take the empty set. Extend the set to $E\subset V$ by including $\{s<-\frac{L}{2}\}$ and excluding $\{s>\frac{L}{2}\}$. Its interior perimeter support $\Sigma$ is compact in $\Int U$. In every dimension the singular set has Hausdorff dimension at most $n-8$, and stability holds for compactly supported variations on $\Sigma_{\reg}$ \cite[Proposition~4.6]{BZ26}.

With $\nu$ the outward unit normal, apply \cite[Lemma~6.8]{Cho21} on $\Sigma_{\reg}$. Since $h'(s)<0$ and $ds(\nu)\le1$, the potential term in that inequality satisfies
$$
R_g+\frac{n}{n-1}h(s)^2+2\,dh(\nu)
\ge \sigma+\frac{n}{n-1}h(s)^2+2h'(s)
=\sigma-\frac{4\pi^2(n-1)}{nL^2}
\ge \sigma-\frac{4\pi^2(n-1)}{n\ell^2}.
$$
This proves the stated estimate.
\end{proof}

\begin{proof}[Proof of \cref{prop:sys-width}]
Suppose the asserted bound fails. Choose $\ell<\width(V,g)$ so that $$\tau=\sigma-4\pi^2\frac{n-1}{n\ell^2}>0, $$and take $E,\Sigma,D$ from \cref{lem:band-separator}.

Choose an SYS class $z$ on $Y$, a detector $u\in H^2(Y; A)$ for $z$, and the torus map $\Theta:Y\to\T^{n-3}$ representing its defining one-classes, as in \cref{prop:detector,prop:carrier}. Write $F=(p,t), G=\Theta\circ p$, and orient $P=\partial_-V$ so that $p_*[P]_\Z=[Y]_\Z$. 

Choose a common regular value $y$ of $G$ on $V$, $\partial V$ and $\Sigma_{\reg}$ outside $G(D)$. Such values have full measure since $\dim_\HH G(D)\le n-8<n-3$.

The fibre $G^{-1}(y)$ is a compact three-manifold transverse to $\Sigma_{\reg}$ and disjoint from $D$. The part of this fibre on the $E$ side, including its boundary, is therefore a compact oriented bordism between $P_y=P\cap G^{-1}(y)$ and $\Sigma_y=\Sigma_{\reg}\cap G^{-1}(y)$. With compatible fibre orientations, $[\Sigma_y]=[P_y]$ in $H_2(V;\Z)$. The projection formula yields
$$
p_*[P_y]=p_*[P]\frown\Theta^*\eta_{\T^{n-3}}=z,\qquad
\langle p^*u,[\Sigma_y]\rangle=\langle u,z\rangle\ne0.
$$

Write $g_\Sigma=g|_{T\Sigma_{\reg}}$. For every $\varphi\in C_c^\infty(\Sigma_{\reg})$, \cref{lem:band-separator} implies
$$
2Q_{\frac{n}{n+1},\tau}^{g_\Sigma,0}(\varphi)=\int_{\Sigma_{\reg}}\bigl(2|d\varphi|_{g_\Sigma}^2+(R_{g_\Sigma}-\tau)\varphi^2\bigr)\,d\HH_g^{n-1}\ge0.
$$
Since $E$ is locally minimizing for the prescribed mean curvature functional near $\Sigma$, the local estimates in the proof of \cref{lem:singular-blowup-potential} apply at $D$. The construction in the proof of \cref{prop:next-metric}, with $m=n$, initial weight zero and loss $\frac{\tau}{2}$, yields a complete metric $g'$ on each component of $\Sigma_{\reg}$ and a smooth weight $f'$ satisfying
$$
g'\ge g_\Sigma,\qquad R_{g'}+2\Delta_{g'}f'-\frac{n-1}{n}|df'|_{g'}^2\ge\frac{\tau}{2}.
$$
For each $a\in D$ there exist $c_a,r_a>0$ such that
$$
g'(x)\ge\frac{c_a}{d_g(x,a)^2}\,g_\Sigma(x)\qquad \forall x\in\Sigma_{\reg}\cap B_g(a,r_a).
$$
If $D=\varnothing$, the same construction takes $g'=g_\Sigma$.

 For the universal covering $\pi_\T:\R^{n-3}\to\T^{n-3}$, set
$$
\widehat V=\{(v,\xi)\in V\times\R^{n-3}:G(v)=\pi_\T(\xi)\},\qquad \pi(v,\xi)=v,\quad \Phi(v,\xi)=\xi.
$$
Both squares below are pullbacks.
$$
\begin{tikzcd}[column sep=large,row sep=large]
\pi^{-1}(\Sigma) \arrow[r,hook] \arrow[d,"\pi"'] & \widehat V \arrow[r,"\Phi"] \arrow[d,"\pi"] & \R\mathrlap{{}^{n-3}} \arrow[d,"\pi_\T"] \\
\Sigma \arrow[r,hook] & V \arrow[r,"G"'] & \T\mathrlap{{}^{n-3}}
\end{tikzcd}
$$
Choose a complete metric $\bar g\ge g$ on $\Int V$ equal to $g$ near $\Sigma$ by taking $\bar g=d_g(\cdot,\partial V)^{-2}g$ in a boundary collar and interpolating. Then
$$
d\Phi_i=\pi^*G^*d\theta_i,\qquad \sup_{\Int\widehat V}|d\Phi_i|_{\pi^*\bar g}\le\sup_V|G^*d\theta_i|_g<\infty.
$$
Fix $\widehat y\in\pi_\T^{-1}(y)$. The projection $\pi$ identifies the oriented fibre $\pi^{-1}(\Sigma)\cap\Phi^{-1}(\widehat y)=\Sigma_y\times\{\widehat y\}$ with $\Sigma_y$. For each component $X$ of $\Int\widehat V$, put $\Gamma_X=X\cap(\Sigma_y\times\{\widehat y\})$ and $u_X=(p\circ\pi|_X)^*u$. Only finitely many components meet this compact fibre, and
$$
\sum_X\langle u_X|_{\Gamma_X},[\Gamma_X]_\Z\rangle=\langle p^*u,[\Sigma_y]_\Z\rangle=\langle u,z\rangle\ne0.
$$
Fix a component $X$ with nonzero summand.

Equip $X$ with $g_X=\pi^*\bar g$. By naturality $u_X$ annihilates $\Sph_2(X)$. Set
$$
X_{n-1}^\circ=X\cap\pi^{-1}(\Sigma_{\reg}),\qquad \overline{X_{n-1}}=X\cap\pi^{-1}(\Sigma),\qquad D_{n-1}=X\cap\pi^{-1}(D).
$$
Since $\bar g=g$ near $\Sigma$ and $\pi$ is a local isometry for the lifted metrics, $\dim_\HH D_{n-1}\le n-8$ and the lifted pair $(\pi^*g',\pi^*f')$ satisfies \ref{item:metric} and \ref{item:inverse}.

Set $F_{n-1}=\Phi|_{\overline{X_{n-1}}}$. For every compact $K\subset\R^{n-3}$, $F_{n-1}^{-1}(K)$ is closed in $\Sigma\times K$ and hence compact. With $V_{n-1}=\R^{n-3}\setminus F_{n-1}(D_{n-1})$ and $N_{n-1}=F_{n-1}^{-1}(V_{n-1})$, we have $\widehat y\in V_{n-1}$ and $F_{n-1}^{-1}(\widehat y)=\Gamma_X$. Thus \cref{prop:fibre degree} yields
$$
\udeg(F_{n-1}|_{N_{n-1}})=\langle u_X|_{\Gamma_X},[\Gamma_X]_\Z\rangle\ne0.
$$
Hence $(X_{n-1}^\circ,\overline{X_{n-1}},D_{n-1},F_{n-1};\pi^*g',\pi^*f',\frac{\tau}{2})$ is an $(n-1)$-dimensional datum, contrary to \cref{prop:datum-obstruction}.

Finally, fix a metric $h$ on $Y$. On $Y\times[-1,1]$, with interval coordinate $t$, take Gromov's warped band metric \cite[pp.~653--654]{Gro18} in the form
$$
g_\varepsilon=dt^2+\left(\frac{\cos(\varepsilon t)}{\cos^2\varepsilon}\right)^{\frac4n}h,\qquad 0<\varepsilon<\frac{\pi}{2}.
$$
Its width is $2$, and the warped-product formula yields
$$
\inf\limits_V R_{g_\varepsilon}=\frac{4(n-1)}n\varepsilon^2+
\left(\frac{\cos^2\varepsilon}{\cos(\varepsilon t)}\right)^{\frac4n}\inf\limits_YR_h\to \frac{(n-1)\pi^2}{n}, \qquad\varepsilon\to\frac{\pi}{2}.
$$
\end{proof}

Following Gromov \cite[Sections~4--6]{Gro18}, define the \emph{SYS-width} of $(M^n,g)$ by
$$
\width_{\mathrm{SYS}}(M,g)=\sup_{\iota:V^n\looparrowright M^n}\width(V,\iota^*g),
$$
where $\iota$ ranges over smooth immersions of compact SYS bands. We call $(M,g)$ \emph{SYS-enlargeable} if $\width_{\mathrm{SYS}}(M,g)=\infty$.

\begin{corollary}\label{cor:sys-enlargeable}
If $(M,g)$ is SYS-enlargeable, then $\inf_M R_g\le0$. In particular, a closed manifold which is SYS-enlargeable for one metric admits no metric of positive scalar curvature.
\end{corollary}
\begin{proof}
If $R_g\ge\sigma>0$, then $R_{\iota^*g}=R_g\circ\iota\ge\sigma$ for every such immersion. Taking the supremum in \cref{prop:sys-width} yields
$$
\width_{\mathrm{SYS}}(M,g)\le2\pi\sqrt{\frac{n-1}{n\sigma}}.
$$
On a closed $M$, any other metric $g_1$ satisfies $g_1\ge c^2g$ for some $c>0$. Hence $\width_{\mathrm{SYS}}(M,g_1)\ge c\,\width_{\mathrm{SYS}}(M,g)$, and the claim follows by compactness.
\end{proof}

\section*{Acknowledgments}
The author would like to express his sincere gratitude to his advisor, Professor Yuguang Shi, for his generous guidance, encouragement, and support throughout this work. The author also thanks Tianze Hao for helpful discussions and suggestions.

\section*{AI Use Statement}
ChatGPT-6 Astra was used to polish and revise the manuscript, help work out the details of the homological comparison via an integral-current filling and some details when constructing the bubble.

\section*{Conflict of Interest}
On behalf of all authors, the corresponding author states that there is no conflict of interest.

\section*{Data Availability}
No datasets were generated or analysed during the current study.

\bibliographystyle{amsalpha-nodash}
\bibliography{reference}
\end{document}